\documentclass{article}%
\usepackage{graphicx}
\usepackage{amsmath}
\usepackage{amscd}
\usepackage{amsfonts}
\usepackage{amssymb}

\newtheorem{theorem}{Theorem}[section]
\newtheorem{lemma}[theorem]{Lemma}
\newtheorem{proposition}[theorem]{Proposition}

\newtheorem{definition}[theorem]{Definition}
\newtheorem{corollary}[theorem]{Corollary}
\newtheorem{remark}[theorem]{Remark}

\newenvironment{proof}[1][Proof]{\textbf{#1.} }{\ \rule{0.5em}{0.5em}}

\newcommand{\be}{\begin{equation}}
\newcommand{\ee}{\end{equation}}
\newcommand{\bes}{\begin{equation*}}
\newcommand{\ees}{\end{equation*}}

\newcommand{\cA}{\mathcal{A}}

\newcommand{\cE}{\mathcal{E}}

\newcommand{\cH}{\mathcal{H}}
\newcommand{\cK}{\mathcal{K}}
\newcommand{\cM}{\mathcal{M}}
\newcommand{\cN}{\mathcal{N}}
\newcommand{\cF}{\mathcal{F}}
\newcommand{\cL}{\mathcal{L}}

\newcommand{\cP}{\mathcal{P}}
\newcommand{\cQ}{\mathcal{Q}}
\newcommand{\cR}{\mathcal{R}}
\newcommand{\cS}{\mathcal{S}}
\newcommand{\tT}{\tilde{T}}
\newcommand{\tV}{\tilde{V}}
\newcommand{\mfp}{\mathfrak{p}}
\newcommand{\mfq}{\mathfrak{q}}

\newcommand{\mfB}{\mathfrak{B}}
\newcommand{\Rp}{\mathbb{R}_+}
\newcommand{\Rpt}{\mathbb{R}_+^2}

\begin{document}

\title{E$_0$-dilation of strongly commuting CP$_0$-semigroups}

\author{Orr Moshe Shalit}
\maketitle
\begin{abstract}
{We prove that every strongly commuting pair of CP$_0$-semigroups
has a minimal E$_0$-dilation. This is achieved in two major steps,
interesting in themselves: {\bf 1:} we show that a strongly
commuting pair of CP$_0$-semigroups can be represented via a two
parameter product system representation; {\bf 2:} we prove that
every fully coisometric product system representation has a fully
coisometric, isometric dilation. In particular, we obtain that
every commuting pair of CP$_0$-semigroups on $B(H)$, $H$ finite dimensional,
has an E$_0$-dilation.}
\end{abstract}

\section{Introduction}
Let $\cM$ be a von Neumann algebra acting on a separable Hilbert space $H$. A \emph{CP$_0$-semigroup on $\cM$} is a family $\Theta = \{\Theta_t\}_{t\geq0}$ of normal, unital, completely positive maps on $\cM$ satisfying the semigroup property
$$\Theta_{s+t}(a) = \Theta_s (\Theta_t(a)) \,\, ,\,\, s,t\geq 0, a\in \cM ,$$
and the continuity condition
$$\lim_{t\rightarrow t_0} \langle\Theta_t(a)h,g\rangle = \langle\Theta_{t_0}(a)h,g\rangle \,\, , \,\, a\in \cM, h,g \in H .$$
A CP$_0$-semigroup is also called a \emph{Quantum Markov Process}, as it may be considered as a noncommutative analog of a Markov process.

A CP$_0$-semigroup is called an \emph{E$_0$-semigroup} if each of
its elements is a $*$-endomorphism. In the past two decades,
E$_0$-semigroups have been extensively studied (for a thorough
introduction, including many references and ``historical" remarks,
see \cite{Arv03}). Although every E$_0$-semigroup is a
CP$_0$-semigroup, non-multiplicative CP$_0$-semigroups are known
to be quite different from E$_0$-semigroups. However, it has been
proved that, in some sense, every CP$_0$-semigroup is ``part" of a
bigger E$_0$-semigroup. To be more precise, we say that a
quadruple $(K,u,\cR,\alpha)$ is an \emph{E$_0$-dilation} of
$\Theta$ if $K$ is a Hilbert space, $u: H \rightarrow K$ is an
isometry, $\cR$ is a von Neumann algebra satisfying $u^* \cR u =
\cM$, and $\alpha$ is an E$_0$-semigroup such that \bes \Theta_t
(u^* b u) = u^* \alpha_t (b) u \,\, , \,\, b \in \cR \ees for all
$t \geq 0$. It has been proved by several authors, using several
different techniques, that every CP$_0$-semigroup has an
E$_0$-dilation (Bhat--Skeide \cite{BS00}, SeLegue \cite{SeLegue},
Muhly--Solel \cite{MS02} and Arveson \cite{Arv03}. We note that
most of the authors have this result also for not necessarily
unital semigroups). This is the precise sense in which we mean
that every CP$_0$-semigroup is ``part" of an E$_0$-semigroup.


If $\cS$ is a topological semigroup, one can define the notions of
CP$_0$ and E$_0$-semigroups over $\cS$. It is then natural to ask
wether every CP$_0$-semigroup $\Theta = \{\Theta_s\}_{s\in\cS}$
over $\cS$ has an E$_0$-dilation. In this paper we make a first
attempt to prove the existence of a minimal E$_0$-dilation for a
CP$_0$-semigroup over $\Rpt := [0,\infty)\times [0,\infty)$. Let us now describe
what we actually achieve.

If $\{R_t\}_{t\geq0}$ and $\{S_t\}_{t\geq0}$ are two
CP$_0$-semigroups that commute (that is, for all $t,s\geq0$, $R_s
S_t = S_t R_s$) then we can define a two parameter
CP$_0$-semigroup $P_{(s,t)} = R_s S_t$. And if we begin with a
CP$_0$-semigroup $\{P_{(t,s)}\}_{(t,s)\in \Rpt}$, then we can
define a commuting pair of semigroups by $R_t = P_{(t,0)}$ and
$S_t = P_{(0,t)}$ (there are some non-trivial continuity issues to
take care of. This will be done below, in Lemma
\ref{lem:continuity1}). The main result of this paper is the
following theorem:

\noindent \,\,

\noindent
{\bf Theorem. }\emph{
Let $\{R_t\}_{t\geq0}$ and $\{S_t\}_{t\geq0}$ be two strongly
commuting CP$_0$-semigroups on a von Neumann algebra $\cM\subseteq
B(H)$, where $H$ is a separable Hilbert space. Then the two
parameter CP$_0$-semigroup $P$ defined by
$$P_{(s,t)} := R_s S_t$$
has a minimal E$_0$-dilation.}

\noindent \,\,

The condition of \emph{strong commutativity} that appears in the
above theorem is a technical one, and it is not yet completely
understood. However, there are many pairs of strongly commuting
CP$_0$-semigroups, and in the appendix we give some sufficient,
and in some cases even necessary, conditions for strong
commutativity. These give rise to many examples of two-parameter semigroups for which the above theorem
applies. In particular, by Proposition \ref{prop:SC_H_finite} below, if $H$ is finite dimensional then
every pair of commuting CP maps on $B(H)$ commute strongly, so every pair of commuting CP$_0$-semigroups on $B(H)$ has a minimal E$_0$-dilation (Corollary \ref{cor:M_n(C)}).

Let us now give an overview of the paper, which should also give some idea of how the above theorem is proved. In what follows, we shall use the notation of the theorem stated above.

After reviewing some preliminary notions and setting the notation in Section \ref{sec:pre}, we explain in Section \ref{sec:MS} the approach of Muhly and Solel to dilation theory (as it appeared in \cite{MS02}). This is the approach that we will be using.

In Section \ref{sec:rep_SC}, after introducing the notion of strong commutativity and proving a few related results, we
construct a
(discrete) product system of $\cM'$-correspondences $X$ over $\Rpt$,
together with a fully coisometric, completely contractive covariant representation $(\sigma,T)$
of $X$ on $H$, such that for all $a \in \cM, (s,t) \in \Rpt$,
\bes
P_{(s,t)} (a) = \tilde{T}_{(s,t)} \left(I_{X(s,t)}
\otimes a \right) \tilde{T}_{(s,t)}^* .
\ees
It is in the construction of the product system $X$ that strong commutativity plays its role.

In Section \ref{sec:iso_dil}, we prove that every fully
coisometric, completely contractive covariant representation of a
product system over $\mathbb{R}^k_+$ (and over some more general
semigroups, as well) can be dilated to an \emph{isometric} and
fully coisometric covariant representation. We do this using the
method of ``representing product system representations as
contractive semigroups on a Hilbert space", which we have
introduced in \cite{Shalit07}.

In Section \ref{sec:Edilation} we show that the isometric dilation
$(\rho,V)$ of the product system representation $(\sigma,T)$
obtained in Section \ref{sec:iso_dil} gives rise to our sought
after E$_0$-dilation in the following way (up to a few
simplifications that we must make here):

Let $K$ be the Hilbert space on which $V$ represents $X$, put $\cR
= \rho(M')'$, and let $u$ be the isometric inclusion $H
\rightarrow K$. The E$_0$-dilation we are looking for is
$(K,u,\cR,\alpha)$, where the semigroup
$\alpha=\{\alpha_s\}_{s\in\Rpt}$ is defined by \bes \alpha_s(b) =
\tilde{V}_s(I \otimes b)\tilde{V}_s^* \,\, , \,\, s\in\Rpt , b\in
\cR. \ees At the end of section \ref{sec:Edilation} we show that
the dilation that we constructed is minimal, and we show that if
$\cM = B(H)$ then $\cR = B(K)$.

In Section \ref{sec:prospects} we close this paper by considering
some possible directions for future research. We give an example
of a finding an E$_0$-dilation to a CP$_0$-semigroup over
$\mathbb{N}\times \mathbb{R}_+$ where strong commutativity does
not occur. We also briefly sketch our program for dilating $k$
strongly commuting, unital CP maps.

\section{Preliminaries}\label{sec:pre}
\subsection{$C^*$/$W^*$-correspondences, their products and their representations}\label{subsec:correspondences}
\begin{definition}
Let $\cA$ be a $C^*$-algebra. A \emph{Hilbert
$C^*$-correspondences over $\cA$} is a (right) Hilbert
$\cA$-module $E$ which carries an adjointable, left action of
$\cA$.
\end{definition}
\begin{definition}
Let $\cM$ be a $W^*$-algebra. A \emph{Hilbert
$W^*$-correspondences over $\cM$} is a self-adjoint Hilbert
$C^*$-correspondence $E$ over $\cM$, such that the map $\cM
\rightarrow \cL(E)$ which gives rise to the left action is normal.
\end{definition}
The following notion of representation of a $C^*$-correspondence was studied extensively in \cite{MS98}, and turned out to be a very useful tool.
\begin{definition}
Let $E$ be a $C^*$-correspondence over $\cA$, and let $H$ be a
Hilbert space. A pair $(\sigma, T)$ is called a \emph{completely
contractive covariant representation} of $E$ on $H$ (or, for
brevity, a \emph{c.c. representation}) if
\begin{enumerate}
    \item $T: E \rightarrow B(H)$ is a completely contractive linear map;
    \item $\sigma : A \rightarrow B(H)$ is a nondegenerate $*$-homomorphism; and
    \item $T(xa) = T(x) \sigma(a)$ and $T(a\cdot x) = \sigma(a) T(x)$ for all $x \in E$ and  all $a\in\cA$.
\end{enumerate}
If $\cA$ is a $W^*$-algebra and $E$ is $W^*$-correspondence then
we also require that $\sigma$ be normal.
\end{definition}
Given a $C^*$-correspondence $E$ and a c.c. representation
$(\sigma,T)$ of $E$ on $H$, one can form the Hilbert space $E
\otimes_\sigma H$, which is defined as the Hausdorff completion of
the algebraic tensor product with respect to the inner product
$$\langle x \otimes h, y \otimes g \rangle = \langle h, \sigma (\langle x,y\rangle) g \rangle .$$
One then defines $\tilde{T} : E \otimes_\sigma H \rightarrow H$ by
$$\tilde{T} (x \otimes h) = T(x)h .$$

As in the theory of contractions on a Hilbert space, there are
certain particularly nice representations which deserve to be
singled out.
\begin{definition}
A c.c. representation $(\sigma, T)$ is called \emph{isometric} if
for all $x, y \in E$,
\begin{equation*}
T(x)^*T(y) = \sigma(\langle x, y \rangle) .
\end{equation*}
(This is the case if and only if $\tilde{T}$ is an isometry.) It
is called \emph{fully coisometric} if $\tilde{T}$ is a coisometry.
\end{definition}

Given two Hilbert $C^*$-correspondences $E$ and $F$ over $\cA$,
the \emph{balanced} (or \emph{inner}) tensor product $E
\otimes_{\cA} F$ is a Hilbert $C^*$-correpondence over $\cA$
defined to be the Hausdorf completion of the algebraic tensor
product with respect to the inner product
$$\langle x \otimes y, w \otimes z \rangle = \langle y , \langle x,w\rangle \cdot z \rangle \,\, , \,\,  x,w\in E, y,z\in F .$$
The left and right actions are defined as $a \cdot (x \otimes y) =
(a\cdot x) \otimes y$ and $(x \otimes y)a = x \otimes (ya)$,
respectively, for all $a\in A, x\in E, y\in F$. We shall usually
omit the subscript $\cA$, writing just $E \otimes F$. When working
in the context of $W^*$-correspondences, that is, if $E$ and $F$
are $W$*-correspondences and $\cA$ is a $W^*$-algebra, then $E
\otimes_{\cA} F$ is understood do be the \emph{self-dual
extension} of the above construction.

Suppose $\cS$ is an abelian cancellative semigroup with identity
$0$ and $p : X \rightarrow \cS$ is a family of
$W^*$-correspondences over $\cA$. Write $X(s)$ for the
correspondence $p^{-1}(s)$ for $s \in \cS$. We say that $X$ is a
(discrete) \emph{product system} over $\cS$ if $X$ is a semigroup,
$p$ is a semigroup homomorphism and, for each $s,t \in \cS
\setminus \{0\}$, the map $X(s) \times X(t) \ni (x,y) \mapsto xy
\in X(s+t)$ extends to an isomorphism $U_{s,t}$ of correspondences
from $X(s) \otimes_{\cA} X(t)$ onto $X(s+t)$. The associativity of
the multiplication means that, for every $s,t,r \in \cS$,
\begin{equation}\label{eq:assoc_prod}
U_{s+t,r} \left(U_{s,t} \otimes I_{X(r)} \right) = U_{s,t+r} \left(I_{X(s)} \otimes U_{t,r} \right).
\end{equation}
We also require that $X(0) = \cA$ and that the multiplications
$X(0) \times X(s) \rightarrow X(s)$ and $X(s) \times X(0)
\rightarrow X(s)$ are given by the left and right actions of $\cA$
and $X(s)$.

\begin{definition}
Let $H$ be a Hilbert space, $\cA$ a $W^*$-algebra and $X$ a
product system of Hilbert $\cA$-correspondences over the semigroup
$\cS$. Assume that $T : X \rightarrow B(H)$, and write $T_s$ for
the restriction of $T$ to $X(s)$, $s \in \cS$, and $\sigma$ for
$T_0$. $T$ (or $(\sigma, T)$) is said to be a \emph{completely
contractive covariant representation} of $X$ if
\begin{enumerate}
    \item For each $s \in \cS$, $(\sigma, T_s)$ is a c.c. representation of $X(s)$; and
    \item $T(xy) = T(x)T(y)$ for all $x, y \in X$.
\end{enumerate}
$T$ is said to be an isometric (fully coisometric) representation if it is an isometric (fully coisometric) representation on every fiber $X(s)$.
\end{definition}
Since we shall not be concerned with any other kind of representation, we shall call a completely contractive covariant representation of a product system simply a \emph{representation}.

\subsection{CP-semigroups and E-dilations}

Let $\cS$ be a unital subsemigroup of $\Rp^k$, and let $\cM$ be a
von Neumann algebra acting on a Hilbert space $H$. A \emph{CP map}
is a completely positive, contractive and normal map on $\cM$. A
\emph{CP-semigroup over $\cS$} is a family $\{\Theta_s\}_{s \in
\cS}$ of CP maps on $\cM$ such that
\begin{enumerate}
    \item For all $s,t \in \cS$
    $$\Theta_s \circ \Theta_t = \Theta_{s + t} \,;$$
    \item $\Theta_0 = {\bf id}_{\cM}$;
    \item For all $h,g\in H$ and all $a\in \cM$, the function
    $$\cS\ni s \mapsto \langle \Theta_s(a) h,g \rangle $$
    is continuous.
\end{enumerate}
A CP-semigroup is called an \emph{E-semigroup} if it consists of $*$-endomorphisms. A CP (E) - semigroup is
called a \emph{CP$_0$ (E$_0$)-semigroup} if all its elements are unital.

\begin{definition}\label{def:dilation}
Let $\cM$ be a von Neumann algebra of operators acting on a
Hilbert space $H$, and let $\Theta = \{\Theta_s\}_{s \in \cS}$ be
a CP-semigroup over the semigroup $\cS$. An \emph{E-dilation of
$\Theta$} is a quadruple $(K,u,\cR,\alpha)$, where $K$ is a
Hilbert space, $u: H \rightarrow K$ is an isometry, $\cR$ is a von
Neumann algebra satisfying $u^* \cR u = \cM$, and $\alpha$ is an
E-semigroup over $\cS$ such that \be\label{eq:CPdef_dil1} \Theta_s
(u^* a u) = u^* \alpha_s (a) u \,\, , \,\, a \in \cR \ee for all
$s \in \cS$.

If $(K,u,\cR,\alpha)$ is a dilation of $\Theta$, then $(\cM,
\Theta)$ is called a \emph{compression} of $(K,u,\cR,\alpha)$.
\end{definition}

Let us review some basic facts regarding E-dilations. Most of the
content of the following paragraphs is spelled out in Chapter 8 of
\cite{Arv03}, for the case where $\cS = \mathbb{R}_+$.

Note that by putting $a = u x u^*$ in (\ref{eq:CPdef_dil1}), for
any $x \in \cM$, one has \be\label{eq:CPdef_dil2} \Theta_s (x) =
u^* \alpha_s (u x u^*) u \,\, , \,\, x \in \cM. \ee

If one identifies $\cM$ with $u \cM u^*$, $H$ with $u H$, and $p$
with $u u^*$, one may give the following equivalent definition,
which we shall use interchangeably with definition
\ref{def:dilation}: \emph{a triple $(p,\cR,\alpha)$ is called a
\emph{dilation} of $\Theta$ if $\cR$ is a von Neumann algebra
containing $\cM$, $\alpha$ is an E-semigroup on $\cR$ and $p$ is a
projection in $\cR$ such that $\cM = p \cR p$ and
$$\Theta_s (p a p) = p \alpha_s (a) p $$
holds for all $s \in \cS, a \in \cR$.}

With this change of notation, we have
$$p\alpha_s (a) p = \Theta_s (pap) = \Theta_s (p^2 a p^2) = p\alpha_s (pap) p, $$
so, taking $a = 1-p$,
$$0 = p \alpha_s (p (1-p) p) p = p \alpha_s (1-p) p .$$
This means that for all $s \in \cS$, $\alpha_s (1-p) \leq 1 - p$.
A projection with this property is called \emph{coinvariant} (note
that if $\alpha$ is an E$_0$-semigroup then $p$ is a coinvariant
projection if and only if it is increasing, i.e., $\alpha_s(p)
\geq p$ for all $s\in \cS$). Equivalently,
\be\label{eq:CPdef_dil3} u u^* \alpha_s(1) = u u^* \alpha_s (u
u^*) \,\, , \,\, s \in \cS . \ee One can also show that
(\ref{eq:CPdef_dil2}) and (\ref{eq:CPdef_dil3}) together imply
(\ref{eq:CPdef_dil1}), and this leads to another equivalent
definition of E-dilation of a CP-semigroup.

Let $\Theta = \{\Theta_s \}_{s \in \cS}$ be a CP-semigroup on a
von Neumann algebra $\cM$, and let $(K,u,\cR,\alpha)$ be an
E-dilation of $\Theta$. Assume that $q \in \cR$ is a projection
satisfying $u u^* \leq q$. Assume furthermore that $q$ is
coninvariant. Then one can show that the maps
$$\beta_s : a \mapsto  q \alpha_s (a) q $$
are the elements of a CP-semigroup on $q \cR q$.

If the maps $\{\beta_s\}$ happen to be multiplicative on $q\cR q$,
then we  say that $q$ is \emph{multiplicative}. In this case,
$(qK,u, q \cR q, \beta)$ is an E-dilation of $\Theta$, which is in
some sense ``smaller" than $(K,u,\cR,\alpha)$.

On the other hand, consider the von Neumann algebra
$$\tilde{\cR} =
W^*\left(\bigcup_{s \in \cS} \alpha_s (u \cM u^*)\right) .$$ This
algebra is clearly invariant under $\alpha$, and it contains $u
\cM u^*$. Thus, restricting $\alpha$ to $\tilde{\cR}$, we obtain a
``smaller" dilation. This discussion leads to the following
definition.
\begin{definition}\label{def:min_dil}
Let $(K,u,\cR,\alpha)$ be an E-dilation of the CP-semigroup
$\Theta$.  $(K,u,\cR,\alpha)$ is said to a \emph{minimal} dilation
if there is no multiplicative, coinvariant projection $1 \neq q
\in \cR$ such that $u u^* \leq q$, and if
\be\label{eq:W^*-generator} \cR = W^*\left(\bigcup_{s \in \cS}
\alpha_s (u \cM u^*)\right) . \ee
\end{definition}

In \cite{Arv03} Arveson defines a minimal dilation slightly
differently:
\begin{definition}\label{def:min_dil_Arv}
Let $(K,u,\cR,\alpha)$ be an E-dilation of the CP-semigroup
$\Theta$.  $(K,u,\cR,\alpha)$ is said to a \emph{minimal} dilation
if the central support of $u u^*$ in $\cR$ is $1$, and if
(\ref{eq:W^*-generator}) holds.
\end{definition}
The two definitions have been shown to be equivalent in the case
where $\Theta$ is a CP$_0$-semigroup over $\Rp$ (\cite{Arv03}, Section 8.9). We now treat the general case.

\begin{proposition}\label{prop:equiv_def_min}
Definition \ref{def:min_dil} holds if \ref{def:min_dil_Arv} does.
\end{proposition}
\begin{proof}
Assume that Definition
\ref{def:min_dil} is violated. If (\ref{eq:W^*-generator}) is
violated, then Definition \ref{def:min_dil_Arv} is, too. So assume
that (\ref{eq:W^*-generator}) holds, and that there is a
multiplicative, coinvariant projection $1 \neq q \in \cR$ such
that $u u^* \leq q$. Denote $\cA = \{\alpha_s (a) : a \in u \cM
u^* , s \in \cS \}$. By a trivial generalization of Proposition
8.9.4 of \cite{Arv03}, $q$ commutes with $\alpha_s(q \cR q)$ for
all $s \in \cS$, so $q$ commutes with $\cA$, thus $q$ commutes
with $W^*(\bigcup_{s \in \cS} \alpha_s (u \cM u^*))$. In
other words, $q$ is central in $\cR$.
\end{proof}

Wether or not the two definitions are equivalent remains an interesting open question.
To prove that they are, it would be enough to show that the central support of $p = uu^*$  in
$W^*\left(\bigcup_{s \in \cS} \alpha_s (u \cM u^*)\right)$ is a coinvariant projection, because the central support is clearly a multiplicative projection. This has been done by Arveson in Proposition 8.3.4, \cite{Arv03}, for the case of a CP$_0$-semigroup over $\cS = \Rp$. Arveson's proof makes use of the order structure of $\Rp$ and cannot be extended to the case $\Rpt$ with which we are concerned in this paper.

\section{Overview of the Muhly--Solel approach to dilation}\label{sec:MS}
In this section we describe the approach of Muhly and Solel to
dilation of CP-semigroups on von Neumann algebras. This approach
was used by Muhly and Solel to dilate CP-semigroups over
$\mathbb{N}$ and $\mathbb{R}_+$ (\cite{MS02}), and later by Solel
for semigroups over $\mathbb{N}^2$ (\cite{S06}). Our program is to
adapt this approach for semigroups over $\cS = \Rpt$.

\subsection{The basic strategy}\label{subsec:strategy}
Let $\Theta$ be a CP-semigroup over the semigroup $\cS$, usually
acting on a von Neumann algebra $\cM$ of operators in $B(H)$. The
dilation is carried out in two main steps. In the first step, a
(discrete) product system of $\cM'$-correspondences $X$ over $\cS$
is constructed, together with a c.c. representation $(\sigma,T)$
of $X$ on $H$, such that for all $a \in \cM, s \in \cS$,
\be\label{eq:rep_rep}
\Theta_s (a) = \tilde{T_s} \left(I_{X(s)}
\otimes a \right) \tilde{T_s}^* ,
\ee
where $T_s$ is the
restriction of $T$ to $X(s)$. In Proposition 2.21, \cite{MS02}, it
is proven that for any c.c. representation $(\sigma, T)$ of a
$W^*$-correspondence $\cE$ over a $W^*$-algebra $\cN$, the mapping
$a \mapsto \tilde{T_s} \left(I_{X(s)} \otimes a \right)
\tilde{T_s}^*$ is a normal, completely positive map on
$\sigma(\cN)'$ (for all $s$). It is also shown that if $T$ is
isometric then this map is multiplicative. Having this in mind,
one sees that a natural way to continue the process of dilation
will be to ``dilate" $(\sigma, T)$ to an isometric c.c.
representation.
\begin{definition}
Let $\cA$ be a $C^*$-algebra, $X$ be a product system of
$\cA$-correspondences over the semigroup $\cS$, and $(\sigma, T)$
a c.c. representation of $X$ on a Hilbert space $H$. An
\emph{isometric dilation} of $(\sigma, T)$ is an isometric
representation $(\rho, V)$ of $X$ on a Hilbert space $K \supseteq
H$, such that
\begin{itemize}
    \item[(i)] $H$ reduces $\rho$ and $\rho(a)\big|_H = P_H \rho(a)\big|_H = \sigma (a)$, for all $a\in\cA$;
    \item[(ii)] for all $s \in \cS, x \in X_s$, one has   $P_H V_s (x)\big|_{K \ominus H} = 0$;
    \item[(iii)] for all $s \in \cS, x \in X_s$, one has  $P_H V_s (x)\big|_{H} = T_s (x)$.
\end{itemize}
\end{definition}
Such a dilation is called \emph{minimal} in case the smallest
subspace of $K$ containing $H$ and invariant under every $V_s
(x)$, $x \in X, s \in \cS$, is all of $K$.

It will be convenient at times to regard an isometric dilation as
a quadruple $(K,u,V,\rho)$, where $(\rho, V)$ are as above and
$u:H \rightarrow K$ is an isometry.

Constructing a minimal isometric dilation $(K,u,V,\rho)$ of the
representation $(\sigma, T)$ appearing in equation
(\ref{eq:rep_rep}) constitutes the second step of the dilation
process. Then one has to show that if $\cR = \rho(\cM')'$, and
$\alpha$ is defined by \bes \alpha_s (a) := \tilde{V_s}
\left(I_{X(s)} \otimes a \right) \tilde{V_s}^* \,\,,\,\, a \in \cR
, \ees then the quadruple $(K,u,\cR,\alpha)$ is an E-dilation for
$(\Theta, \cM)$. In \cite{MS98}, \cite{MS02} and \cite{S06}, it is
proved that any c.c. representation of a product system over
$\mathbb{N}$, $\mathbb{R}_+$ or $\mathbb{N}^2$ (the latter two,
$X$ is assumed to be a product system of $W^*$-correspondence, and
$\sigma$ is assumed to be normal), has a minimal isometric
dilation. Moreover, it is shown that if $X$ is a product system of
$W^*$-correspondences and $\sigma$ is assumed to be normal then
$\rho$ is also normal. When the product system is over
$\mathbb{N}$ or $\mathbb{R}_+$, the minimal isometric dilation is
also unique. From these results, the authors deduce the existence
of an E-dilation of a CP-semigroup $\Theta$ acting on a von
Neumann $\cM$. When $\Theta$ is a CP-semigroup over $\cS = \mathbb{R}_+$ and $H$
is seperable, then $\alpha$ is shown to be an
$E$-semigroup that is a minimal dilation.

\subsection{Description of the construction of the product system and representation for one parameter semigroups}\label{subsec:des_MS}
In this subsection we give a detailed description of Muhly and
Solel's construction of the product system and c.c. representation
associated with a one-parameter CP-semigroup (\cite{MS02}). We
shall use this construction in section \ref{sec:rep_SC}. We note that
the original construction in \cite{MS02} was carried out for CP$_0$-semigroups,
it works just as well for CP-semigroups.

Let $\Theta = \{\Theta_t \}_{t \geq 0}$ be a CP-semigroup
acting on a von Neumann algebra $\cM$ of operators in $B(H)$
(we will not really use any assumptions regarding the continuity with respect to
$t$). Let $\mfB(t)$ denote the collection of partitions of the
closed unit interval $[0,t]$, ordered by refinement. For $\mfp \in
\mfB(t)$, we define a Hilbert space $H_{\mfp,t}$ by
$$H_{\mfp,t} : = \cM \otimes_{\Theta_{t_1}} \cM \otimes_{\Theta_{t_2 - t_1}} \cM \otimes \cdots \otimes_{\Theta_{t-t_{n-1}}} H ,$$
where $\mfp = \{0 = t_0 < t_1 < t_2 < \cdots < t_n = t \}$, and the RHS of the above equation is the Hausdorff completion of the algebraic tensor product $\cM \otimes \cM \otimes \cdots \otimes H$ with respect to the inner product
\begin{align*}\langle T_1 \otimes \cdots \otimes T_n \otimes h, & S_1 \otimes \cdots \otimes S_n \otimes k \rangle = \\
& \langle h, \Theta_{t-t_{n-1}} (T_n^* \Theta_{t_{n-1}-t_{n-2}} (T_{n-1}^* \cdots \Theta_{t_1} (T_1^* S_1) \cdots S_{n-1}) S_n) k \rangle .
\end{align*}
$H_{\mfp,t}$ is a left $\cM$-module via the action $S \cdot (T_1 \otimes \cdots \otimes T_n \otimes h) = ST_1 \otimes \cdots \otimes T_n \otimes h$. We now define the intertwining spaces
$$\cL_{\cM} (H,H_{\mfp,t}) = \{X\in B(H,H_{\mfp,t}): \forall S \in \cM . X S = S \cdot X \} .$$
The inner product
$$\langle X_1, X_2 \rangle := X_1^* X_2 ,$$
for $X_i \in \cL_{\cM} (H,H_{\mfp,t})$, together with the right and left actions
$$(XR)h: = X(Rh) ,$$
and
$$(RX)h := (I \otimes \cdots \otimes I \otimes R) Xh ,$$
for $R \in \cM', X \in \cL_{\cM} (H,H_{\mfp,t})$, make $\cL_{\cM} (H,H_{\mfp,t})$ into a $W^*$-correspondence over $\cM'$.

The Hilbert spaces $H_{\mfp,t}$ and $W^*$-correspondences $\cL_{\cM} (H,H_{\mfp,t})$ form inductive systems as follows. Let $\mfp,\mfp' \in \mfB(t)$, $\mfp \leq \mfp'$. In the particular case where $\mfp = \{0 = t_0 < \cdots < t_k < t_{k+1}< \cdots < t_n = t\}$ and $\mfp' = \{0 = t_0 < \cdots < t_k < \tau < t_{k+1}< \cdots < t_n = t\}$, we can define a Hilbert space isometry $v_0 : H_{\mfp,t} \rightarrow H_{\mfp',t}$ by
\begin{align*}
v_0 (T_1 \otimes \cdots \otimes T_{k+1} \otimes T_{k+2} \otimes & \cdots \otimes T_n \otimes h) = \\
& T_1 \otimes \cdots \otimes T_{k+1} \otimes I \otimes T_{k+2} \otimes \cdots \otimes T_n \otimes h.
\end{align*}
This map gives rise to an isometry of $W^*$-correspondences $v : \cL_{\cM} (H,H_{\mfp,t}) \rightarrow \cL_{\cM} (H,H_{\mfp',t})$ by $v(X) = v_0 \circ X$.

Now, if $\mfp \leq \mfp'$ are any partitions in $\mfB(t)$, then we can define $v_{0,\mfp,\mfp'} : H_{\mfp,t} \rightarrow H_{\mfp',t}$ and $v_{\mfp,\mfp'} : \cL_{\cM} (H,H_{\mfp,t}) \rightarrow \cL_{\cM} (H,H_{\mfp',t})$ by composing a finite number of maps such as $v_0$ and $v$ constructed in the previous paragraph, and we geéçt legitimate arrow maps. Now one can form two different direct limits:
$$H_t := \underrightarrow{\lim }(H_{\mfp,t},v_{0,\mfp,\mfp'}) $$
and
$$E(t) := \underrightarrow{\lim }(\cL_{\cM} (H,H_{\mfp,t}),v_{\mfp,\mfp'}) .$$
The inductive limit also supplies us with embeddings of the blocks $v_{0,\mfp,\infty} : H_{\mfp,t} \rightarrow H_{t}$ and $v_{\mfp,\infty} : \cL_{\cM} (H,H_{\mfp,t}) \rightarrow E(t)$. One can also define interwining spaces $\cL_{\cM} (H,H_{t})$, each of which has the structure of an $\cM'$-correspondence, and these spaces are isomorphic as $W^*$-correspondences to the spaces $E(t)$. $\{E(t)\}_{t \geq 0}$ is the product sytem of $\cM'$-correspondences that we are looking for. We have yet to describe the the c.c. representation $(\sigma, T)$ that will ``represent" $\Theta$ as in
equation (\ref{eq:rep_rep}) (with $X(s)$ replaced by $E(s)$).

The sought after representation is the so called ``identity representation", which we now describe. First, we set $\sigma = T_0 = {\bf id}_{\cM'}$. Next, let $t > 0$. For $\mfp = \{0 = t_0 < \cdots < t_n = t\}$, the formula
$$\iota_\mfp (h) = I \otimes \cdots \otimes I \otimes h$$
defines an isometry $\iota_\mfp : H \rightarrow H_{\mfp,t}$, with adjoint given by the formula
$$\iota_{\mfp}^*(X_1 \otimes \cdots \otimes X_n \otimes h) = \Theta_{t-t_{n-1}}(\Theta_{t_{n-1}-t_{n-2}}(\cdots (\Theta_{t_1}(X_1)X_2) \cdots X_{n-1})X_n)h .$$
For $\mfp'$ a refinement of $\mfp$, one computes $\iota_{\mfp}^* = \iota_{\mfp'}^* \circ v_{0,\mfp,\mfp'}$. This induces a unique map $\iota_t^* : H_t \rightarrow H$ that satisfies $\iota_t^* \circ v_{0,\mfp,\infty} = \iota_{\mfp}^*$. The c.c. representation $T_t$ on $E(t)$ is given by
$$T_t (X) = \iota_t^* \circ X ,$$
where we have identified $E(t)$ with $\cL_{\cM} (H,H_{t})$.

\section{Representing strongly commuting CP$_0$-semigroups}\label{sec:rep_SC}
In this section and in the next two we prove our main result:
\emph{every pair of strongly commuting CP$_0$-semigroups has an
E$_0$-dilation}. As we mentioned in the previous section, our
program is to prove this result using the Muhly-Solel approach,
which consists of two main steps. In this section we concentrate
on the first step: the representation of a pair of strongly
commuting CP-semigroups using a product system representation
via a formula such as equation (\ref{eq:rep_rep}) above. This will
be done in the third subsection, whereas the first and second
subsections will be devoted to the notion of strong commutativity
and its implications.

Throughout this and the two following sections, $\cM$ will be a von Neumann algebra
acting on a Hilbert space $H$. There is a natural correspondence
between two parameter semigroups of maps and pairs of commuting
one parameter semigroups. Indeed, if $\{R_t\}_{t\geq0}$ and
$\{S_t\}_{t\geq0}$ are two semigroups that commute (that is, for
all $t,s\geq0$, $R_s S_t = S_t R_s$) then we can define a two
parameter semigroup $P_{(s,t)} = R_s S_t$. And if we begin with a
semigroup $\{P_{(t,s)}\}_{(t,s)\in \Rpt}$, then we can define a
commuting pair of semigroups by $R_t = P_{(t,0)}$ and $S_t =
P_{(0,t)}$. It is not trivial that $P$ is continuous (in the relevant sense)
 if and only if
$R$ and $S$ are -- it follows from the fact that $(s,X) \mapsto
R_s(X)$ is jointly contiuous in the weak topology (we shall make
this argument precise in Lemma \ref{lem:continuity1}). From now on
we fix the notation in the preceding paragraph, and we shall use
either $\{P_{(t,s)}\}_{(t,s)\in \Rpt}$ or the pair
$\{R_t\}_{t\geq0}$ and $\{S_t\}_{t\geq0}$ to denote a fixed
two-parameter CP-semigroup. Note also that if
$\{\alpha_t\}_{t\geq0}$ and $\{\beta_t\}_{t\geq0}$ are
\emph{commuting} E-dilations of $\{R_t\}_{t\geq0}$ and
$\{S_t\}_{t\geq0}$ acting on the same von Neumann algebra, then
$\{\alpha_t \beta_s\}_{t,s\geq0}$ is an E-dilation of
$\{P_{(t,s)}\}_{(t,s)\in \Rpt}$, and vice versa.

\subsection{Strongly commuting CP maps}
Let $\Theta$ and $\Phi$ be CP maps on $\cM$. We define the Hilbert
space $\cM \otimes_\Phi \cM \otimes_\Theta H$ to be the Hausdorff
completion of the algebraic tensor product $\cM
\otimes_\textrm{alg} \cM \otimes_\textrm{alg} H$ with respect to
the inner product
$$\langle a \otimes b \otimes h, c \otimes d \otimes k \rangle = \langle h, \Theta(b^* \Phi(a^* c) d) k \rangle .$$
\begin{definition}\label{def:SC}
Let $\Theta$ and $\Phi$ be CP maps on $\cM$. We say that they
\emph{commute strongly} if there is a unitary $u: \cM \otimes_\Phi
\cM \otimes_\Theta H \rightarrow \cM \otimes_\Theta \cM
\otimes_\Phi H$ such that:
\begin{itemize}
\item[(i)] $u(a \otimes_\Phi I \otimes_\Theta h) = a \otimes_\Theta I \otimes_\Phi h$ for all $a \in \cM$ and $h \in H$.
\item[(ii)] $u(ca\otimes_\Phi b \otimes_\Theta h) = (c \otimes I_M \otimes I_H)u(a\otimes_\Phi b \otimes_\Theta h)$ for $a,b,c \in \cM$ and $h \in H$.
\item[(iii)] $u(a\otimes_\Phi b \otimes_\Theta dh) = (I_M \otimes I_M \otimes d)u(a\otimes_\Phi b \otimes_\Theta h)$ for $a,b \in \cM$, $d \in \cM'$ and $h \in H$.
\end{itemize}
\end{definition}
The notion of strong commutation was introduced by Solel in
\cite{S06}. Note that if two CP maps commute strongly, then they
commute. The converse is false (for concrete examples see Subsections
\ref{subsec:NSC} and \ref{subsec:SC_comm}). In
the appendix we shall give many examples of strongly
commuting pairs of CP maps, and for some von Neumann algebras we shall
give a complete characterization of strong commutativity. For the time being
let us just state the fact that if $H$ is a finite dimensional Hilbert space,
then any two commuting CP maps on $B(H)$ strongly commute (see Subsection \ref{subsec:H_finite}).
The ``true" significance of strong
commutation comes from a bijection between pairs of strongly
commuting CP maps and product systems over $\mathbb{N}^2$ with
c.c. representations (\cite{S06}, Propositions 5.6 and 5.7, and
the discussion between them). It is this bijection that enables one
to characterize all pairs of strongly commuting CP maps on $B(H)$
(\cite{S06}, Proposition 5.8).

In the next section we will work with the spaces $\cM
\otimes_{P_1} \cM  \cdots  \cM \otimes_{P_n} H$, where $P_1,
\ldots, P_n$ are CP maps. These spaces are defined in a way
analogous to the way that the spaces $\cM \otimes_\Theta \cM
\otimes_\Phi H$ were defined in the beginning of this section. The
following results are important for dealing with such spaces.

\begin{lemma}\label{lem:SC1}
Assume that $P_{n-1}$ and $P_n$ commute strongly. Then there exists a unitary
$$v : \cM \otimes_{P_1} \cM \otimes_{P_2} \cdots \otimes_{P_{n-1}} \cM \otimes_{P_n} H
\rightarrow
\cM \otimes_{P_1} \cM \otimes_{P_2} \cdots \otimes_{P_n} \cM \otimes_{P_{n-1}} H
$$
such that
\begin{enumerate}
    \item $v(I \otimes_{P_1} \cdots \otimes_{P_{n-1}} I \otimes_{P_n} h) = I \otimes_{P_1} \cdots \otimes_{P_n} I \otimes_{P_{n-1}} h$, for all $h \in H$,
    \item For all $X \in \cM$,
    $$v \circ (X \otimes I \cdots I \otimes I)
    = (X \otimes I \cdots  I \otimes I)\circ v ,$$
    \item For all $X \in \cM'$,
    $$v \circ (I \otimes I \cdots I \otimes X)
    = (I \otimes I \cdots  I \otimes X)\circ v .$$
\end{enumerate}
\end{lemma}
\begin{proof}
Let $u : \cM \otimes_{P_{n-1}} \cM \otimes_{P_n} H \rightarrow \cM
\otimes_{P_n} \cM \otimes_{P_{n-1}} H$ be the unitary that makes
$P_{n-1}$ and $P_n$ commute strongly. Define
$$v = I_E \otimes u ,$$
where $E$ denotes the $W^*$-correspondence (over $\cM$) $\cM
\otimes_{P_1} \cM \otimes_{P_2} \cdots \otimes_{P_{n-3}} \cM$
equipped with the inner product
$$\langle a_1 \otimes \cdots \otimes a_{n-3}, b_1 \otimes \cdots \otimes b_{n-3} \rangle
= P_{n-3}\left(a_{n-3}^*\cdots P_1(a_1^* b_1) \cdots
b_{n-3}\right) .$$

The fact that $v$ commutes with $\cM \otimes I \otimes \cdots
\otimes I$ and $I \otimes I \cdots I \otimes \cM'$ and satisfies
the three conditions listed above are clear from the definition
and from the properties of $u$. The fact that $u$ is surjective
implies that $v$ is, too. It is left to show that $v$ is an
isometry (and this will also show that it is well defined). Let
$\sum a_i \otimes_{P_{n-2}} b_i \otimes_{P_{n-1}} c_i
\otimes_{P_n} h_i$ be an element of $E \otimes_{P_{n-2}} \cM
\otimes_{P_{n-1}} \cM \otimes_{P_n} H$.
\begin{align*}
& \|v (\sum a_i \otimes_{P_{n-2}} b_i \otimes_{P_{n-1}} c_i \otimes_{P_n} h_i) \|^2
 = \\
& = \langle \sum a_i \otimes_{P_{n-2}} u(b_i \otimes_{P_{n-1}} c_i \otimes_{P_n} h_i), \sum a_j \otimes_{P_{n-2}} u(b_j \otimes_{P_{n-1}} c_j \otimes_{P_n} h_j) \rangle
 = \\
& = \sum_{i,j} \langle u(b_i \otimes_{P_{n-1}} c_i \otimes_{P_n} h_i), P_{n-2}\left(\langle a_i, a_j \rangle\right) u(b_j \otimes_{P_{n-1}} c_j \otimes_{P_n} h_j) \rangle
 = (*) \\
& = \sum_{i,j} \langle u(b_i \otimes_{P_{n-1}} c_i \otimes_{P_n} h_i), u\left(P_{n-2}\left(\langle a_i, a_j \rangle\right) b_j \otimes_{P_{n-1}} c_j \otimes_{P_n} h_j\right) \rangle
 = (**) \\
& = \sum_{i,j} \langle b_i \otimes_{P_{n-1}} c_i \otimes_{P_n} h_i, P_{n-2}(\langle a_i, a_j \rangle) b_j \otimes_{P_{n-1}} c_j \otimes_{P_n} h_j \rangle
 = \\
& =\|\sum a_i \otimes_{P_{n-2}} b_i \otimes_{P_{n-1}} c_i \otimes_{P_n} h_i \|^2
\end{align*}
the equality marked by (*) follows from the fact that $u$
interwines the actions of $\cM$ on $\cM \otimes_{P_{n-1}} \cM
\otimes_{P_{n}} H$ and $\cM \otimes_{P_{n}} \cM \otimes_{P_{n-1}}
H$, and the one marked by (**) is true because $u$ is unitary.
\end{proof}

\begin{lemma}\label{lem:SC2}
Assume that $P$ and $Q$ are strongly commuting CP maps on $\cM$.
Then there exists an isomorphism $v = v_{P,Q}$ of $\cM$-correspondences
$$v : \cM \otimes_P \cM \otimes_Q \cM \rightarrow \cM \otimes_Q \cM \otimes_P \cM $$
such that
$$v(I \otimes_P I \otimes_Q I) = I \otimes_Q I \otimes_P I .$$
\end{lemma}

\begin{remark}
\emph{In the sequel, given strongly commuting CP maps $P$ and $Q$,
it will be convenient to refer to the $\cM$-module isometry
$v_{P,Q}$ as the \emph{associated map}.}
\end{remark}
\begin{proof}
For any two CP maps $\Theta, \Phi$ let $W_{\Theta,\Phi}$ be the
Hibert space isomorphism
$$W_{\Theta, \Phi} : \cM \otimes_\Theta \cM \otimes_\Phi \cM \otimes_I H \rightarrow  \cM \otimes_\Theta \cM \otimes_\Phi H$$
given by $W_{\Theta, \Phi} (a\otimes_\Theta b \otimes_\Phi c
\otimes_I h) = a \otimes_\Theta b \otimes_\Phi ch$. By a
straightforward computation $W_{\Theta, \Phi}^*$ is given by
$W_{\Theta, \Phi}^* (a \otimes_\Theta b \otimes_\Phi h) = a
\otimes_\Theta b \otimes_\Phi I \otimes_I h$, and by even shorter
computations $W_{\Theta, \Phi}W_{\Theta, \Phi}^*$ and $W_{\Theta,
\Phi}^*W_{\Theta, \Phi}$ are identity maps. For all $a,b,c,x \in
\cM$ and all $y \in \cM'$ we have
\begin{align*}\label{eq:W_properties}
W_{\Theta, \Phi}(xa\otimes_\Theta b \otimes_\Phi c \otimes_I yh) & =
xa\otimes_\Theta b \otimes_\Phi cyh \\ & =
x a\otimes_\Theta b \otimes_\Phi ych \\ & =
(x \otimes I \otimes y)W_{\Theta, \Phi}(a\otimes_\Theta b \otimes_\Phi c \otimes_I h) .
\end{align*}
From this, it also follows that
$$W_{\Theta, \Phi}^*(x \otimes I \otimes y) = (x \otimes I \otimes I \otimes y)W_{\Theta, \Phi}^* \quad \,\, (x \in \cM, y \in \cM') .$$

We now define a map $T : \cM \otimes_P \cM \otimes_Q \cM \otimes_I
H \rightarrow \cM \otimes_Q \cM \otimes_P \cM \otimes_I H$ by
$$T = W_{Q,P}^* \circ u \circ W_{P,Q} ,$$
where $u$ is the map that makes $P$ and $Q$ commute strongly. As a
product of such maps, $T$ is a unitary interwining the left
actions of $\cM$ and $\cM'$. The $v$ that we are looking for is a
map $v : \cM \otimes_P \cM \otimes_Q \cM \rightarrow \cM \otimes_Q
\cM \otimes_P \cM$ that satisfies $T = v \otimes I_H$. We will
find this $v$ using a standard technique exploiting the self
duality of $\cM \otimes_Q \cM \otimes_P \cM$.

For any $x \in \cM \otimes_Q \cM \otimes_P \cM$ we define a map $
L_x : H \rightarrow \cM \otimes_Q \cM \otimes_P \cM \otimes_I H$
by
$$L_x (h) = x \otimes h \,\, , \,\,(h \in H) .$$
The adjoint is given on simple tensors by $L_x^*(y \otimes h) = \langle x, y\rangle h$.

Now, if there is a $v$ such that $T = v \otimes I_H$, then for all $z\in \cM \otimes_P \cM \otimes_Q \cM$ and $x\in \cM \otimes_Q \cM \otimes_P \cM$ we must have
$$
\langle x, v(z) \rangle h = L_x^* (v(z) \otimes h)  = L_x^* T(z \otimes h) .
$$
This leads us to define, fixing $z \in \cM \otimes_P \cM \otimes_Q
\cM$, a mapping $\varphi$ from $\cM \otimes_Q \cM \otimes_P \cM$
into $\cM$:
$$\varphi (x) h := L_x^*T (z \otimes h) .$$
We now prove that $x \mapsto \varphi(x)^*$ is a bounded, $\cM$-module mapping into $\cM$.

{\bf into $\cM$:} For all $x \in \cM \otimes_Q \cM \otimes_P \cM$,
$\varphi(x)$ is linear. $\|L_x^* T (z \otimes h) \| \leq \|L_x^*
\|\|T \|\|z \|\|h \| $, so $\varphi(x) \in B(H)$. So
$\varphi(x)^*$ exists and is also a bounded, linear operator on
$H$. Now take $d \in \cM'$. Then
$$\varphi(x)dh = L_x^*T (z \otimes dh) = L_x^*T (I \otimes d)(z \otimes h) = L_x^* (I \otimes d) T (z \otimes h) = d \varphi(x)h$$
($L_x^*$ interwines $\cM'$ from its very definition) whence $\varphi(x) \in \cM'' = \cM$. Thus, $\varphi(x)^* \in \cM$.

{\bf $\cM$-module mapping:} This is because for all $x,y \in \cM
\otimes_Q \cM \otimes_P \cM$ and all $a \in \cM$ $L_{x+y} = L_x +
L_y$ and $L_{ax} = a L_x$ (and also $L_{xa} = L_x a$).

{\bf Bounded mapping:} From the inequalities $\|L_x^* T (z \otimes
h) \| \leq \|L_x^* \|\|T \|\|z \|\|h \| $ and $\|L_x^*\| \leq
\|x\|$ it follows that $\|\varphi(x)^*\|=\|\varphi(x)\| \leq \|z\|
\|x\|$.

It now follows from the self-duality of $\cM \otimes_Q \cM
\otimes_P \cM$ that for all $ z \in \cM \otimes_P \cM \otimes_Q
\cM$ there exists a $v(z) \in \cM \otimes_Q \cM \otimes_P \cM$
such that \be\label{eq:existv} \langle x,v(z)\rangle h = L_x^* T(z
\otimes h) \ee for all $x \in \cM \otimes_Q \cM \otimes_P \cM, h
\in H$. It is easy to see from (\ref{eq:existv}) that $v(z)$ is a
right $\cM$-module mapping. (\ref{eq:existv}) can be re-written as
$$L_x^* (v(z) \otimes h) = L_x^* T(z \otimes h) ,$$
and, since this holds for all $x$, this means that $(v(z) \otimes
h) = T(z \otimes h)$ (because $\cap_x {\rm Ker}(L_x^*) =
\left(\vee_x {\rm Im}(L_x)\right)^\perp = \{0\}$),
or, in other words, $v
\otimes I = T$. This last equality implies that $v$ is unitary,
and that it has all the properties required. For example, if
$a,b,c,X \in \cM$ and $h \in H$, then
\begin{align*}
v(Xa \otimes b \otimes c) \otimes h & = T(Xa \otimes b \otimes c \otimes h) \\
& = (X \otimes I \otimes I \otimes I)T(a \otimes b \otimes c \otimes h) \\
& = (X \otimes I \otimes I \otimes I)(v(a \otimes b \otimes c)\otimes h) \\
& = \big((X \otimes I \otimes I)(v(a \otimes b \otimes c)\big) \otimes h .
\end{align*}
Putting $v_1 = v(Xa \otimes b \otimes c)$ and $v_2=(X \otimes I \otimes I)(v(a \otimes b \otimes c)$ we have that for all $h \in H$
$$0 = \|v_1 \otimes h - v_2 \otimes h \|^2= \|(v_1 - v_2)\otimes h\|^2 = \langle h, \langle v_1 - v_2,v_1 - v_2 \rangle h \rangle ,$$
which implies that $\langle v_1 - v_2,v_1 - v_2 \rangle = 0$, or $v(Xa \otimes b \otimes c) =(X \otimes I \otimes I)(v(a \otimes b \otimes c)$.
\end{proof}

\begin{remark}\label{rem:SC}
\emph{The converse of Lemma \ref{lem:SC2} is also true: if there is an isometry of $\cM$-correspondences $v : \cM \otimes_P \cM \otimes_Q \cM \rightarrow \cM \otimes_Q \cM \otimes_P \cM$ such that $v(I \otimes I \otimes I) = I \otimes I \otimes I$ then $P$ and $Q$ strongly commute. Indeed, to obtain $u : \cM \otimes_P \cM \otimes_Q H \rightarrow \cM \otimes_Q \cM \otimes_P H$ with the desired properties, we simply reverse the construction above. That is, we define $T = v \otimes I$, and}
$$u = W_{Q,P} \circ T \circ W_{P,Q}^* .$$
\end{remark}

\begin{lemma}\label{lem:SC3}
Assume that $P_j$ and $P_{j+1}$ commute strongly, for some $j\leq n-2$. Then there exists a unitary
$$u : \cM \otimes_{P_1} \cdots \otimes_{P_j} \cM \otimes_{P_{j+1}} \cdots \cM \otimes_{P_n} H
\rightarrow
\cM \otimes_{P_1} \cdots \otimes_{P_{j+1}} \cM \otimes_{P_{j}} \cdots \cM \otimes_{P_n} H
$$
such that
\begin{enumerate}
    \item $u(I \otimes_{P_1} \cdots I \otimes_{P_j} I \otimes_{P_{j+1}}I \cdots I \otimes_{P_n} h) = I \otimes_{P_1} \cdots I \otimes_{P_{j+1}} I \otimes_{P_j}I \cdots I \otimes_{P_n} h$,
    \item For all $X \in \cM$,
    $$u \circ (X \otimes I \cdots I \otimes I)
    = (X \otimes I \cdots  I \otimes I)\circ u ,$$
    \item For all $X \in \cM'$,
    $$u \circ (I \otimes I \cdots I \otimes X)
    = (I \otimes I \cdots  I \otimes X)\circ u .$$
\end{enumerate}
\end{lemma}
\begin{proof}
Let $v : \cM \otimes_{P_{j}} \cM \otimes_{P_{j+1}} \cM \rightarrow
\cM \otimes_{P_{j+1}} \cM \otimes_{P_{j}} \cM$ be the unitary that
is described in lemma \ref{lem:SC2}. Introduce the notation
$$E = \cM \otimes_{P_{1}} \cdots \otimes_{P_{j-2}} \cM$$
(understood to be $\mathbb{C}$ if $j=1$ and $\cM$ if $j=2$) and
$$F = \cM \otimes_{P_{j+3}} \cdots \cM \otimes_{P_{n}} H$$
(understood to be $H$ if $j=n-2$). Define
$$u : E \otimes_{P_{j-1}} \cM \otimes_{P_{j}} \cM \otimes_{P_{j+1}} \cM \otimes_{P_{j+2}} F
\rightarrow
E \otimes_{P_{j-1}} \cM \otimes_{P_{j+1}} \cM \otimes_{P_{j}} \cM \otimes_{P_{j+2}} F$$
by
$$u := I_E \otimes v \otimes I_F .$$
$u$ is a well-defined, unitary mapping, possesing the properties asserted.
\end{proof}

Putting together Lemmas \ref{lem:SC1}, \ref{lem:SC2} and
\ref{lem:SC3}, we obtain the following
\begin{proposition}\label{prop:SC}
Let $R_1, R_2, \ldots R_m$, and $S_1, S_2, \ldots, S_n$ be CP maps
such that for all $1 \leq i \leq m$, $1 \leq j \leq n$, $R_i$
commutes strongly with $S_j$. Then there exists a unitary
$$v : \cM \otimes_{R_1} \cdots \otimes_{R_m} \cM \otimes_{S_1} \cdots \otimes_{S_n} H
\rightarrow
\cM \otimes_{S_1} \cdots \otimes_{S_n} \cM \otimes_{R_1} \cdots \otimes_{R_m} H$$
such that
\begin{enumerate}
    \item $v(I \otimes_{R_1} I \cdots  I \otimes_{S_n} h) = I \otimes_{S_1} I \cdots I \otimes_{R_{m}} h$, for all $h \in H$,
    \item For all $X \in \cM$,
    $$v \circ (X \otimes I \cdots I \otimes I)
    = (X \otimes I \cdots  I \otimes I)\circ v ,$$
    \item For all $X \in \cM'$,
    $$v \circ (I \otimes I \cdots I \otimes X)
    = (I \otimes I \cdots  I \otimes X)\circ v .$$
\end{enumerate}
\end{proposition}
The existence of $v$ as above is clear: simply apply the
isomorphisms from the previous lemmas one by one. One might think
that applying these isomorphisms in different orders might lead to
different $v$'s. In the next subsection we will see, however, that
the order of application does not influence the total outcome (cf.
Proposition \ref{prop:SCSG}).

\subsection{Strongly commuting CP-semigroups}

\begin{definition}\label{def:SCSG}
Two semigroups of CP maps $\{R_t\}_{t\geq0}$ and
$\{S_t\}_{t\geq0}$ are said to \emph{commute strongly} if for all
$(s,t) \in \Rpt$ the CP maps $R_s$ and $S_t$ commute strongly.
\end{definition}

In the appendix we have collected a few examples of strongly
commuting CP-semigroups, and we give some necessary and sufficient
conditions for strong commutativity in special cases. From this point on
$R$ and $S$ will denote two strongly commuting CP-semigroups.

\begin{proposition}\label{prop:SCSG}
If the CP-semigroups $\{R_t\}_{t\geq0}$ and $\{S_t\}_{t\geq0}$
commute strongly, then, for all $(s,t),(s',t') \in \Rpt$, the
associated maps
$$v_{R_s,S_t}: \cM \otimes_{R_s} \cM \otimes_{S_t} \cM \rightarrow \cM \otimes_{S_t} \cM \otimes_{R_s} \cM ,$$
and
$$v_{R_{s'},S_{t'}}: \cM \otimes_{R_{s'}} \cM \otimes_{S_{t'}} \cM \rightarrow \cM \otimes_{S_{t'}} \cM \otimes_{R_{s'}} \cM ,$$
(see lemma \ref{lem:SC2}) satisfy the following identity :
\be\label{eq:SCSG} (I \otimes I \otimes
v_{R_{s'},S_{t'}})(v_{R_s,S_t} \otimes I \otimes I) = (v_{R_s,S_t}
\otimes I \otimes I)(I \otimes I \otimes v_{R_{s'},S_{t'}}). \ee
\end{proposition}

\begin{proof}
Let $a,b,c,d,e \in \cM$. Assume that $v_{R_s,S_t} (a \otimes_{R_s}
b \otimes_{S_t} c) = \sum_{i=1}^m A_i \otimes_{S_t} B_i
\otimes_{R_s} C_i$, and that $v_{R_{s'},S_{t'}} (I
\otimes_{R_{s'}} d \otimes_{S_{t'}} e) = \sum_{j=1}^n \gamma_i
\otimes_{S_{t'}} \delta_j \otimes_{R_{s'}} \epsilon_j$. Operating
on $a \otimes_{R_s} b \otimes_{S_t} c \otimes_{R_{s'}} d
\otimes_{S_{t'}} e$ with the operator on the LHS of equation
(\ref{eq:SCSG}), we obtain
\begin{align*}
& (I \otimes I \otimes v_{R_{s'},S_{t'}})(v_{R_s,S_t} \otimes I \otimes I)(a \otimes b \otimes c \otimes d \otimes e) = \\
& (I \otimes I \otimes v_{R_{s'},S_{t'}}) \sum_{i=1}^m A_i \otimes B_i \otimes C_i \otimes d \otimes e = (*)\\
& \sum_{i=1}^m A_i \otimes B_i \otimes C_i \cdot v_{R_{s'},S_{t'}} (I \otimes d \otimes e) = \\
& \sum_{i=1}^m \sum_{j=1}^n A_i \otimes B_i \otimes C_i \gamma_j \otimes \delta_j \otimes \epsilon_j ,
\end{align*}
where the equality marked by (*) is justified because
$v_{R_{s'},S_{t'}}$ is a left $\cM$-module map. Operating on $a
\otimes_{R_s} b \otimes_{S_t} c \otimes_{R_{s'}} d
\otimes_{S_{t'}} e$ with the operator on the RHS of equation
(\ref{eq:SCSG}), we obtain
\begin{align*}
& (v_{R_s,S_t} \otimes I \otimes I)(I \otimes I \otimes v_{R_{s'},S_{t'}}) (a \otimes b \otimes c \otimes d \otimes e) = (*) \\
& (v_{R_s,S_t} \otimes I \otimes I)(a \otimes b \otimes c \cdot v_{R_{s'},S_{t'}} (I \otimes d \otimes e)) = \\
& \sum_{j=1}^n (v_{R_s,S_t} \otimes I \otimes I) (a \otimes b \otimes c \gamma_j \otimes \delta_j \otimes \epsilon_j) = \\
& \sum_{j=1}^n v_{R_s,S_t} (a \otimes b \otimes c \gamma_j) \otimes \delta_j \otimes \epsilon_j = (**)\\
& \sum_{j=1}^n v_{R_s,S_t} (a \otimes b \otimes c ) \cdot \gamma_j \otimes \delta_j \otimes \epsilon_j = \\
& \sum_{i=1}^m \sum_{j=1}^n A_i \otimes B_i \otimes C_i \gamma_j \otimes \delta_j \otimes \epsilon_j ,
\end{align*}
where the equality marked by (*) is justified because
$v_{R_{s'},S_{t'}}$ is a left $\cM$-module map, and the one marked
by (**) is OK because $v_{R_{s},S_{t}}$ is a right $\cM$-module
map. So equation (\ref{eq:SCSG}) holds for all $s,s',t,t' \geq$,
and this proof is complete.
\end{proof}

\subsection{Representing a pair of strongly commuting CP$_0$-semigroups
 via the identity representation - the strongly commuting case}\label{subsec:repvia}

Recall the notation that we fixed in this chapter: $\cM$ is a von
Neumann algebra acting on $H$, $\{R_t\}_{t\geq0}$ and
$\{S_t\}_{t\geq0}$ are two strongly commuting CP-semigroups on
$\cM$, and $P_{(s,t)} := R_s S_t$. Let $\{E(t)\}_{t\geq0}$,
$\{F(t)\}_{t\geq0}$ denote the product systems (of
$W^*$-correspondences over $\cM '$) associated with
$\{R_t\}_{t\geq0}$ and $\{S_t\}_{t\geq0}$, respectively, and let
$T^E$, $T^F$ be the corresponding identity representations (as
described in Subsection \ref{subsec:des_MS}). For $s,t \geq 0$, we
denote by $\theta_{s,t}^E$ and $\theta_{s,t}^F$ the isometries
$$\theta_{s,t}^E : E(s)\otimes_{\cM'}E(t) \rightarrow E(s+t) ,$$
and
$$\theta_{s,t}^F : F(s)\otimes_{\cM'}F(t) \rightarrow F(s+t) .$$

\begin{proposition}\label{prop:tech}
For all $s,t \geq 0$ there is an isomorphism of
$W^*$-correspondences \be \varphi_{s,t}:E(s) \otimes_{\cM'} F(t)
\rightarrow F(t) \otimes_{\cM'} E(s). \ee The isomorphisms
$\{\varphi_{s,t}\}_{s,t\geq0}$, together with the identity
represetations $T^E$, $T^F$, satisfy
the ``commutation" relation: \be\label{eq:commutation_relation}
\tT_s^E (I_{E(s)} \otimes \tT_t^F) = \tT_t^F (I_{F(t)} \otimes
\tT_s^E) \circ (\varphi_{s,t} \otimes I_H) \quad , t,s \geq0 . \ee
\end{proposition}
\begin{proof}
We shall adopt the notation used in subsection \ref{subsec:des_MS}
(with a few changes), and follow the proof of proposition 5.6 in
\cite{S06}. Fix $s,t \geq 0$. Let $\mathfrak{p} = \{0 = s_0 < s_1
< \ldots < s_m = s\}$ be a partition of $[0,s]$. We define
$$H_{\mathfrak{p}}^R = \cM \otimes_{R_{s_1}} \cM \otimes_{R_{s_2-s_1}} \cdots \cM \otimes_{R_{s_m-s_{m-1}}} H$$
and we define (for a partition $\mathfrak{q}$) $H_{\mathfrak{q}}^S$ in a similar manner. If $q = \{0 = t_0 < t_1 < \ldots < t_n = t\}$, we also define
$$H_{\mathfrak{p},\mathfrak{q}}^{R,S} = \cM \otimes_{R_{s_1}} \cdots \otimes_{R_{s_m-s_{m-1}}} \cM \otimes_{S_1}
\cdots \otimes_{S_{t_n-t_{n-1}}} H .$$
$H_{\mathfrak{q},\mathfrak{p}}^{S,R}$ is defined similarly. We can
go on to define
$H_{\mathfrak{q},\mathfrak{p},\mathfrak{p}'}^{S,R,S}$,
$H_{\mathfrak{q},\mathfrak{p},\mathfrak{q}',\mathfrak{p}'}^{S,R,S,R}$,
etc.

Recall that $E(s)$ is the direct limit of the directed system
$(\cL_\cM (H,H_{\mathfrak{p}}^R), v_{\mathfrak{p,p'}})$.
Similarily, we shall write $(\cL_\cM (H,H_{\mathfrak{q}}^S),
u_{\mathfrak{q,q'}})$ for the directed system that has $F(t)$ as
its limit. We write $v_{\mathfrak{p},\infty}$,
$u_{\mathfrak{q},\infty}$ for the limit isometric embeddings.

We proceed to construct an isomorphism
$$\varphi_{s,t} : E(s)\otimes F(t) \rightarrow F(t) \otimes E(s)$$
that has the desired property. Let $\mathfrak{p}  = \{0 = s_0 <
s_1 < \ldots < s_m = s\}$ and $\mathfrak{q}  = \{0 = t_0 < t_1 <
\ldots < t_n = t\}$ be partitions of $[0,s]$ and $[0,t]$,
respectively. Denote by $\Gamma_{\mfp,\mfq}$ The map from $\cL_\cM
(H,H_{\mathfrak{p}}^R) \otimes \cL_\cM (H,H_{\mathfrak{q}}^S)$
into $\cL_\cM (H,H_{\mathfrak{q},\mathfrak{p}}^{S,R})$ given by $X
\otimes Y \mapsto (I \otimes I \cdots I \otimes X)Y$. As explained
in lemma 3.2 of \cite{MS02}, $\Gamma_{\mfp,\mfq}$ is an
isomorphism. We define $\Gamma_{\mfq,\mfp}$ to be the
corresponding map from $\cL_\cM (H,H_{\mathfrak{q}}^S) \otimes
\cL_\cM (H,H_{\mathfrak{p}}^R)$ into $\cL_\cM
(H,H_{\mathfrak{p},\mathfrak{q}}^{R,S})$. Let $u :
H_{\mathfrak{q},\mathfrak{p}}^{S,R} \rightarrow
H_{\mathfrak{p},\mathfrak{q}}^{R,S}$ be the isomorphism from
corollary \ref{prop:SC}, and define $\Psi : \cL_\cM
(H,H_{\mathfrak{q},\mathfrak{p}}^{S,R}) \rightarrow \cL_\cM
(H,H_{\mathfrak{p},\mathfrak{q}}^{R,S})$ by $\Psi(Z) = u \circ Z$.
The argument from proposition 5.6 from \cite{S06} can be repeated
here to show that $\Psi$ is an isomorphism of
$W^*$-correspondences. Define $t_{\mathfrak{p},\mathfrak{q}} :
\cL_\cM (H,H_{\mathfrak{p}}^R) \otimes \cL_\cM
(H,H_{\mathfrak{q}}^S) \rightarrow \cL_\cM (H,H_{\mathfrak{q}}^S)
\otimes \cL_\cM (H,H_{\mathfrak{p}}^R)$ by
$$t_{\mathfrak{p},\mathfrak{q}} = \Gamma_{\mfq,\mfp}^{-1}\circ \Psi \circ \Gamma_{\mfp,\mfq} .$$
Define maps $W_1 : H \rightarrow H_{\mathfrak{p}}^R$ and $W_2 : H
\rightarrow H_{\mathfrak{q}}^S$ by $W_1 h = I \otimes_{R_1} \cdots
I \otimes_{R_{s_m-s_{m-1}}}  h$ and $W_2 h = I \otimes_{S_1}
\cdots I \otimes_{S_{t_n-t_{n-1}}} h$. Also, let $U_1 :
H_{\mathfrak{p}}^R \rightarrow
H_{\mathfrak{q},\mathfrak{p}}^{S,R}$ and $U_2 : H_{\mathfrak{q}}^S
\rightarrow H_{\mathfrak{p},\mathfrak{q}}^{R,S}$ be the maps $U_1
\xi = I \otimes_{S_1} I \cdots I \otimes_{S_{t_n-t_{n-1}}} \xi$
and $U_2 \eta = I \otimes_{R_1} I \cdots I
\otimes_{R_{s_m-s_{m-1}}} \eta$. Just as in \cite{S06}, we have
that
\begin{equation}\label{eq:WU}
W_1^* U_1^* = W_2^* U_2^* u ,
\end{equation}
and that, for $X \in \cL_\cM (H,H_{\mathfrak{p}}^R)$, we have
$U_1^* (I \otimes \cdots I \otimes X) = X W_2^*$. Now, for $X \in
\cL_\cM (H,H_{\mathfrak{p}}^R)$ and $Y \in \cL_\cM
(H,H_{\mathfrak{q}}^S)$,
\begin{equation}\label{eq:UGamma}
U_1^* \Gamma_{\mfp,\mfq} (X \otimes Y) = U_1^* (I \otimes I \cdots
I \otimes X)Y = X W_2^* Y.
\end{equation}
If we define $(T_\mathfrak{p}^R, id)$ \footnote{Watch out - we
have here a little problem with notation - this resembles $T_t^E,
T_t^F$ that we defined above.} to be the identity representation
of $\cL_\cM (H,H_{\mathfrak{p}}^R)$, and $(T_\mathfrak{q}^S, id)$
to be the identity representation of $\cL_\cM
(H,H_{\mathfrak{q}}^S)$, (see the closing paragraph in subsection
\ref{subsec:des_MS}), then (\ref{eq:UGamma}) implies that, for $h
\in H$, \be\label{eq:REP} W_1^* U_1^* (\Gamma_{\mfp,\mfq}(X
\otimes Y))h = T_\mathfrak{p}^R (X) T_\mathfrak{q}^S (Y)h =
\tilde{T}_\mathfrak{p}^R (I \otimes \tilde{T}_\mathfrak{q}^S)
(X\otimes Y\otimes h). \ee On the other hand, using (\ref{eq:WU})
and an analog of (\ref{eq:REP}),
\begin{align*}
W_1^* U_1^* (\Gamma_{\mfp,\mfq}(X \otimes Y))h & =
W_1^* U_1^* (\Psi^{-1} \Gamma_{\mfq,\mfp} \circ t_{\mathfrak{p},\mathfrak{q}} (X \otimes Y))h \\
& = W_1^* U_1^* u^*(\Gamma_{\mfq,\mfp} \circ t_{\mathfrak{p},\mathfrak{q}} (X \otimes Y))h \\
& = W_2^* U_2^* (\Gamma_{\mfq,\mfp} \circ t_{\mathfrak{p},\mathfrak{q}} (X \otimes Y))h \\
& = \tilde{T}_\mathfrak{q}^S (I \otimes \tilde{T}_\mathfrak{p}^R)(t_{\mathfrak{p},\mathfrak{q}} (X \otimes Y) \otimes h) .
\end{align*}

Let us summarize what we have accumulated up to this point. For
fixed $s,t \geq 0$, and any two partitions $\mathfrak{p}$,
$\mathfrak{q}$ of $[0,s]$ and $[0,t]$, respectively, we have a
Hilbert space isomorphism
$$t_{\mathfrak{p},\mathfrak{q}} : \cL_\cM (H,H_{\mathfrak{p}}^R) \otimes \cL_\cM (H,H_{\mathfrak{q}}^S) \rightarrow \cL_\cM (H,H_{\mathfrak{q}}^S) \otimes \cL_\cM (H,H_{\mathfrak{p}}^R)$$
satisfying \be\label{eq:com} \tilde{T}_\mathfrak{p}^R (I \otimes
\tilde{T}_\mathfrak{q}^S) = \tilde{T}_\mathfrak{q}^S (I \otimes
\tilde{T}_\mathfrak{p}^R)(t_{\mathfrak{p},\mathfrak{q}} \otimes
I_H) . \ee These maps induce an isomorphism
$t_{\mathfrak{p},\infty}: \cL_\cM (H,H_{\mathfrak{p}}^R) \otimes
F(t) \rightarrow F(t) \otimes \cL_\cM (H,H_{\mathfrak{p}}^R)$ that
satisfies \be\label{eq:tpinfty} t_{\mathfrak{p},\infty} (I \otimes
u_{\mathfrak{q},\infty}) = (u_{\mathfrak{q},\infty} \otimes I)
t_{\mathfrak{p},\mathfrak{q}} . \ee Plugging (\ref{eq:tpinfty}) in
(\ref{eq:com}) we obtain
$$
\tilde{T}_\mathfrak{p}^R (I \otimes \tilde{T}_\mathfrak{q}^S) = \tilde{T}_\mathfrak{q}^S (u_{\mathfrak{q},\infty}^* \otimes \tilde{T}_\mathfrak{p}^R)(t_{\mathfrak{p},\infty} \otimes I_H) (I \otimes u_{\mathfrak{q},\infty} \otimes I_H) .
$$
The discussion before theorem 3.9 in \cite{MS02} imply that
$\tilde{T}_t^F (u_{\mathfrak{q},\infty} \otimes I) =
\tilde{T}_{\mathfrak{q}}^S$, or, letting $p_\mathfrak{q}$ denote
the projection in $F(t)$ onto $u_{\mathfrak{q},\infty} (\cL_\cM
(H,H_{\mathfrak{q}}^S))$,
$$\tilde{T}_t^F (p_{\mathfrak{q}} \otimes I) = \tilde{T}_{\mathfrak{q}}^S  (u_{\mathfrak{q},\infty} ^* \otimes I_H).$$
The last two equations sum up to
$$
\tilde{T}_\mathfrak{p}^R (I \otimes \tilde{T}_t^F)(I \otimes p_{\mathfrak{q}} \otimes I_H) = \tilde{T}_t^F (p_{\mathfrak{q}} \otimes \tilde{T}_\mathfrak{p}^R)(t_{\mathfrak{p},\infty} \otimes I_H) (I \otimes p_{\mathfrak{q}} \otimes I_H) ,
$$
which implies, in the limit,
$$
\tilde{T}_\mathfrak{p}^R (I \otimes \tilde{T}_t^F) = \tilde{T}_t^F (I_{F(t)} \otimes \tilde{T}_\mathfrak{p}^R)(t_{\mathfrak{p},\infty} \otimes I_H)  .
$$
Repeating this ``limiting process" in the argument $\mathfrak{p}$,
we obtain a map $t_{\infty,\infty} : E(s) \otimes F(t) \rightarrow
F(t) \otimes E(s)$, which we re-label as $\varphi_{s,t}$, that
satisfies (\ref{eq:commutation_relation}). The above procedure can
be done for all $s,t \geq 0$, giving isomorphisms
$\{\varphi_{s,t}\}$ satisfying the commutation relation
(\ref{eq:commutation_relation}).
\end{proof}

Our aim now is to construct a product system $X$ over $\Rpt$ and a
c.c. representation $T$ of $X$ that will lead to a representation
of $\{P_{(s,t)}\}_{(s,t)\in \Rpt}$ as in equation
(\ref{eq:rep_rep}). Proposition \ref{prop:tech} is a key
ingredient in the proof that the representation that we define
below gives rise to such a representation. But before going into
that we need to carefully construct the product system $X$.

We define
$$X(s,t) := E(s) \otimes F(t) ,$$
and
$$\theta_{(s,t),(s',t')}: X(s,t) \otimes X(s',t') \rightarrow X(s+s',t+t') ,$$
by
$$\theta_{(s,t),(s',t')} =  (\theta_{s,s'}^E \otimes \theta_{t,t'}^F)\circ (I \otimes \varphi_{s',t}^{-1} \otimes I) .$$
To show that $\{ X(s,t) \}_{t,s\geq0}$ is a product system, we
shall need to show that ``the $\theta$'s make the tensor product
into an associative multiplication", or simply: \be\label{eq:ass}
\theta_{(s,t),(s'+s'',t'+t'')} \circ (I \otimes
\theta_{(s',t'),(s'',t'')}) = \theta_{(s+s',t+t'),(s'',t'')} \circ
(\theta_{(s,t),(s',t')} \otimes I) , \ee for
$s,s',s'',t,t',t''\geq 0$.

\begin{proposition}\label{prop:assoc}
$X = \{ X(s,t) \}_{t,s\geq0}$ is a product system. That is,
equation (\ref{eq:ass}) holds.
\end{proposition}
\begin{proof}
The proof is nothing but a straightforward and tedious computation, using Proposition
\ref{prop:SCSG}.

Let $s,s',s'',t,t',t'' \geq 0$, and let $\mfp, \mfp', \mfp'',
\mfq, \mfq', \mfq''$ be partitions of the corresponding intervals.
It is enough to show that the maps on both sides of equation
(\ref{eq:ass}) give the same result when applied to an element of
the form
$$\zeta = X \otimes Y \otimes X' \otimes Y' \otimes X'' \otimes Y'',$$
where $X \in \cL_\cM (H, H^R_\mfp)$, $Y \in \cL_\cM (H,
H^S_\mfq)$, etc. Let us operate first on $\zeta$ with the RHS of
(\ref{eq:ass}).

Now,
$$\theta_{(s,t),(s',t')} (X \otimes Y \otimes X' \otimes Y') =
\left( \theta_{\mfp,\mfp'}^E \otimes \theta_{\mfq,\mfq'}^F \right)
\left(X \otimes t_{\mfp',\mfq}^{-1}(Y \otimes X') \otimes Y'
\right) ,$$
where $\theta_{\mfp,\mfp'}^E$ is the restriction of
$\theta_{s,s'}^E$ to $\cL_\cM (H,H_\mfp^R) \otimes \cL_\cM
(H,H_{\mfp'}^R)$, $\theta_{\mfq,\mfq'}^F$ is defined similarly,
and $t_{\mfp',\mfq}$ is the map defined in Proposition
\ref{prop:tech}. Looking at the definition of $t_{\mfp',\mfq}$, we see
that $t_{\mfp',\mfq}^{-1} (Y \otimes X') =
\Gamma_{\mfp',\mfq}^{-1}\left(U_{\mfp' \leftrightarrow \mfq} \circ
(I\otimes Y) X' \right)$. Here $U_{\mfp' \leftrightarrow \mfq}$
denotes the unitary $H_{\mfp',\mfq}^{R,S} \rightarrow
H_{\mfq,\mfp'}^{S,R}$ given by Proposition \ref{prop:SC}. Assume
that
$$U_{\mfp' \leftrightarrow \mfq} \circ
(I\otimes Y) X' = \sum_{i}(I \otimes x_i)y_i .$$ Then
$$\Gamma_{\mfp',\mfq}^{-1}\left(U_{\mfp' \leftrightarrow \mfq} \circ
(I\otimes Y) X' \right) = \sum_i x_i \otimes y_i ,$$ therefore,
$$\theta_{(s,t),(s',t')} (X \otimes Y \otimes X' \otimes Y') =
\sum_i (I \otimes X)x_i \otimes (I \otimes y_i)Y'.$$
So,
\begin{align*}
& \theta_{(s+s',t+t'),(s'',t'')} \circ (\theta_{(s,t),(s',t')}
\otimes I) \zeta = \\
& \sum_i \left(\theta_{\mfp' \vee \mfp + s', \mfp''}^E \otimes
\theta_{\mfq' \vee \mfq + t', \mfq''}^F \right) \Big[(I \otimes
X)x_i \otimes \Gamma_{\mfp'', \mfq' \vee q + t'}^{-1} \left(
U_{\mfp'' \leftrightarrow \mfq' \vee q + t'} \circ (I \otimes (I \otimes
y_i)Y')X'' \right) \otimes Y'' \Big] .
\end{align*}
Repeated application of Proposition \ref{prop:SCSG} shows that,
and this is a crucial point, $U_{\mfp'' \leftrightarrow \mfq' \vee q + t'} = (I
\otimes U_{\mfp'' \leftrightarrow \mfq}) (U_{\mfp'' \leftrightarrow \mfq'} \otimes I)$. Thus
$$U_{\mfp'' \leftrightarrow \mfq' \vee q + t'} \circ (I \otimes (I \otimes
y_i)Y')X'' = (I \otimes U_{\mfp'' \leftrightarrow \mfq})(I \otimes I \otimes y_i)
\left(U_{\mfp'' \leftrightarrow \mfq'} (I \otimes Y')X'' \right).$$
Write $U_{\mfp'' \leftrightarrow \mfq'} \circ (I \otimes Y')X''$ as $\sum_j (I \otimes a_j)b_j$.
Then we have
\begin{align*}
U_{\mfp'' \leftrightarrow \mfq' \vee q + t'} \circ (I \otimes (I \otimes
y_i)Y')X''
& = \sum_j (I \otimes U_{\mfp'' \leftrightarrow \mfq})(I \otimes I \otimes y_i)
(I \otimes a_j)b_j \\
& = \sum_j \left(I \otimes \big[U_{\mfp'' \leftrightarrow \mfq} \circ (I \otimes y_i)a_j \big]
\right) b_j .
\end{align*}
We now write $U_{\mfp'' \leftrightarrow \mfq} \circ (I \otimes y_i)a_j$ as $\sum_k (I \otimes A_{i,j,k})B_{i,j,k}$.
With this notation, we get
\begin{align*}
& \theta_{(s+s',t+t'),(s'',t'')} \circ (\theta_{(s,t),(s',t')}
\otimes I) \zeta = \\
& \sum_{i,j,k} \left((I \otimes I \otimes X)(I \otimes x_i)A_{i,j,k} \right)
\otimes \left((I \otimes I \otimes B_{i,j,k})(I \otimes b_j)Y'' \right) .
\end{align*}

Now let us operate first on $\zeta$ with the LHS of
(\ref{eq:ass}), repeating all the steps that we have made above:
\begin{align*}
\theta_{(s',t'),(s'',t'')} (X' \otimes Y' \otimes X'' \otimes Y'')
& = \left( \theta_{\mfp',\mfp''}^E \otimes \theta_{\mfq',\mfq''}^F \right)
\left(X' \otimes t_{\mfp'',\mfq'}^{-1}(Y' \otimes X'') \otimes Y''
\right) \\
& = \sum_j (I \otimes X')a_j \otimes (I \otimes b_j)Y'',
\end{align*}
thus,
\begin{align*}
& \theta_{(s,t),(s'+s'',t'+t'')} \circ (I \otimes \theta_{(s',t'),(s'',t'')}) \zeta = \\
& \sum_j \left(\theta_{\mfp,\mfp'' \vee \mfp' + s''}^E \otimes
\theta_{\mfq,\mfq'' \vee \mfq' + t''}^F \right)
\Big[ X  \otimes \Gamma_{\mfp'' \vee \mfp' + s'', \mfq}^{-1} \left(
U_{\mfp'' \vee \mfp' + s'' \leftrightarrow \mfq} \circ (I \otimes (I \otimes
Y)X')a_j \right) \otimes (I \otimes b_j) Y'' \Big] .
\end{align*}
As above, we factor $U_{\mfp'' \vee \mfp' + s'' \leftrightarrow \mfq}$ as
$(U_{\mfp'' \leftrightarrow \mfq} \otimes I) (I \otimes U_{\mfp' \leftrightarrow \mfq})$, to obtain

\begin{align*}
 U_{\mfp'' \vee \mfp' + s'' \leftrightarrow \mfq} \circ (I \otimes (I \otimes
Y)X')a_j  & =
\sum_i (U_{\mfp'' \leftrightarrow \mfq} \otimes I) \circ (I \otimes (I \otimes x_i)y_i)a_j \\
& = \sum_i (I \otimes I \otimes x_i) (U_{\mfp'' \leftrightarrow \mfq} \otimes I) \circ (I \otimes y_i)a_j \\
& = \sum_{i,k} (I \otimes I \otimes x_i) (I \otimes A_{i,j,k})B_{i,j,k} .
\end{align*}
So we get
\begin{align*}
& \theta_{(s,t),(s'+s'',t'+t'')} \circ (I \otimes \theta_{(s',t'),(s'',t'')}) \zeta = \\
& \sum_{i,j,k} \left((I \otimes I \otimes X)(I \otimes x_i)A_{i,j,k} \right)
\otimes \left((I \otimes I \otimes B_{i,j,k})(I \otimes b_j)Y'' \right) ,
\end{align*}
and this is exactly the same expression as we obtained for $\theta_{(s+s',t+t'),(s'',t'')} (\theta_{(s,t),(s',t')} \otimes I) \zeta$.
\end{proof}

\begin{theorem}\label{thm:rep_SC}
There exists a two parameter product system of
$\cM'$-correspondences $X$, and a completely contractive,
covariant representation $T$ of $X$ into
$B(H)$, such that for all $(s,t) \in \Rpt$ and all $a \in \cM$,
the following identity holds: \be\label{eq:rep_SC}
\tT_{(s,t)}(I_{X(s,t)} \otimes a)\tT_{(s,t)}^* = P_{(s,t)}(a) .
\ee
\end{theorem}
\begin{proof}
As above, define
$$X(s,t) := E(s) \otimes F(t) .$$
By Proposition \ref{prop:assoc}, $X$ is a product system. For $s,t
\geq 0$, $\xi \in E(s)$ and $\eta \in F(t)$, we define a representation $T$ of $X$ by
$$T_{(s,t)}(\xi \otimes \eta) := T_s^E (\xi) T_t^F (\eta) .$$
It is clear that for fixed $s,t \geq 0$, $T_{(s,t)}$, together
with $\sigma = {\bf id}_{\cM'}$, extends to a covariant
representation of $X(s,t)$ on $H$. In addition, \be\label{eq:coi}
\tT_{(s,t)} = \tT_s^E (I_{E(s)}\otimes \tT_t^F)  , \ee so
$\|\tT_{(s,t)}\| \leq 1$. By lemma 3.5 in \cite{MS98}, $T_{(s,t)}$
is completely contractive. Also, if $P$ is unital, so are $R$
and $S$, thus $T^E$ and $T^F$ are fully coisometric, whence $T$ is
fully coisometric. We turn to show that for $x_1 \in X(s_1,t_1),
x_2 \in X(s_2,t_2)$, $T_{(s_1 +s_2,t_1 + t_2)} (x_1 \otimes x_2) =
T_{(s_1,t_1)}(x_1) T_{(s_2,t_2)}(x_2)$.

Let $\xi_i \in E(s_i), \eta_i \in F(t_i), \,\, i=1,2$. Put $\Phi =
I_{E(s_1)} \otimes \varphi_{s_2,t_1} \otimes I_{F(t_2)}$. Treating
the maps $\theta_{s_1,s_2}^E , \theta_{t_1,t_2}^F$ as identity
maps, we have that $\Phi : X(s_1 +s_2, t_1 + t_2) \rightarrow
X(s_1,t_1) \otimes X(s_2,t_2)$. We need to show that
$$T_{(s_1 + s_2,t_1 + t_2)}\left(\Phi^{-1}(\xi_1 \otimes \eta_1 \otimes \xi_2 \otimes \eta_2) \right)
= T_{(s_1,t_1)}(\xi_1 \otimes \eta_1) T_{(s_2,t_2)}(\xi_2 \otimes
\eta_2) .$$ But for this it suffices to show that \footnote{when
writing this down, the problems arise only in the ``middle".}
$$T_{(s,t)} \left( \varphi_{s,t}^{-1}(\eta \otimes \xi) \right) = T_{(0,t)}(\eta) T_{(s,0)}(\xi) \quad ,\quad \xi \in E(s), \eta \in F(t) .$$
Let $h \in H$. Now, on the one hand, recalling
(\ref{eq:commutation_relation}), we have
$$\tT_{(s,0)}(I_{E(s)} \otimes \tT_{(0,t)}) (\varphi_{s,t}^{-1}(\eta \otimes \xi) \otimes h) = \tT_{(0,t)}(I_{F(t)} \otimes \tT_{(s,0)})(\eta \otimes \xi \otimes h ) = T_{(0,t)}(\xi) T_{(s,0)}(\eta) h .$$
On the other hand, writing $\sum \xi_i \otimes \eta_i$ for
$\varphi_{s,t}^{-1}(\eta \otimes \xi)$, we have
\begin{align*}
\tT_{(s,0)}(I_{E(s)} \otimes \tT_{(0,t)}) (\varphi_{s,t}^{-1}(\eta
\otimes \xi) \otimes h)
& = \sum \tT_{(s,0)}(\xi_i \otimes T_{(0,t)}(\eta_i) h) \\
& = \sum T_{(s,0)}(\xi_i) T_{(0,t)}(\eta_i) h \\
& = T_{(s,t)}(\sum \xi_i \otimes \eta_i) h\\
& = T_{(s,t)}(\varphi_{s,t}^{-1}(\eta \otimes \xi))h
\end{align*}
so we conclude that $T_{(0,t)}(\xi) T_{(s,0)}(\eta) =
T_{(s,t)}(\varphi_{s,t}^{-1}(\eta \otimes \xi))$, as required.

Finally, using theorem 3.9 in \cite{MS02}, we easily compute for
$a \in \cM$:
\begin{align*}
\tT_{(s,t)}(I_{X(s,t)} \otimes a)\tT_{(s,t)}^* & =
\tT_{(s,0)}(I_{E(s)} \otimes \tT_{(0,t)})(I_{E(s)} \otimes
I_{F(t)} \otimes a)
(I_{E(s)} \otimes \tT_{(0,t)}^*)\tT_{(s,0)}^* \\
& = \tT_{(s,0)}(I_{E(s)} \otimes S_t (a) ) \tT_{(s,0)}^* \\
& = R_s (S_t (a)) = P_{(s,t)} (a) .
\end{align*}
This concludes the proof.
\end{proof}

\section{Isometric dilation of a fully coisometric product system representation}\label{sec:iso_dil}

In the previous section, given a von Neumann algebra $\cM
\subseteq B(H)$ and two strongly commuting CP$_0$-semigroups on
$\cM$, we constructed a product system $X$ of
$\cM'$-correspondences over $\Rpt$ and a product system
representation $(\sigma,T)$ of $X$ on $H$ such that for all $(s,t)
\in \Rpt$ and all $a \in \cM$, \bes \tT_{(s,t)}(I_{X(s,t)} \otimes
a)\tT_{(s,t)}^* = P_{(s,t)}(a) . \ees In other words, we have
completed the first step in our program for dilation. In this
section we shall carry out the second step: we shall construct a
fully coisometric, isometric dilation $(\rho,V)$ of $(\sigma,T)$
on some Hilbert space $K\supseteq H$. In the next section we will
show that the family of maps given by
$$\alpha_{(s,t)}(b) := \tV_{(s,t)}(I_{X(s,t)} \otimes b)\tV_{(s,t)}^* $$
for all $b \in \cR := \rho(\cM')'$ is the E$_0$-dilation that we are looking for.

In fact, we are going to prove a little more than we need: we
shall prove that every fully coisometric representation of a
product system over a (certain kind of) subsemigroup of $\Rp^k$ has an isometric
dilation (see subsection \ref{subsec:coisoisodil}). This result
will be proved by ``representing the representation as a
contractive semigroup on a Hilbert space" (see subsection
\ref{subsec:rep}), a method that we introduced in \cite{Shalit07}.
Since in this paper we shall be ultimately interested in applying
this result for the product system and representation given in
Theorem \ref{thm:rep_SC}, we will not make the construction or
statement in the most general possible way, in hope of making the
presentation as smooth as possible. For example, one does not have
to assume that neither the product system nor the representation
is unital, but we shall make these assumptions, as they hold for
the output of Theorem \ref{thm:rep_SC}. Also, the reader will note
that our construction makes sense for more general semigroups than
those we shall consider.

\subsection{Representing product system representations as contractive semigroups on a Hilbert space}\label{subsec:rep}
Let $\cS$ be a subsemigroup of $\Rp^k$ ($k$ can be taken to be some
infinite cardinal number, but we shall assume $k \in \mathbb{N}$ to keep things simple). Let $\cA$ be a unital $C^*$-algebra,
and let $X$ be a discrete product system of unital
$C^*$-correspondences over $\cS$ \footnote{An $\cA$-correspondence
is said to be \emph{unital} if the left action of $\cA$ is
unital. Note that if $\cA$ is unital, then the right action of $\cA$
on every $\cA$-correspondence is unital.}. Let $(\sigma,T)$ be a completely contractive covariant
representation of $X$ on the Hilbert space $H$, and assume that
$\sigma$ is unital. Our unital assumptions imply that $\cA \otimes
H = X(0) \otimes H \cong H$ via the identification $a \otimes h
\leftrightarrow \sigma(a)h$. This identification will be made repeatedly
below.

Define $\cH_0$ to be the space of all finitely supported functions
$f$ on $\cS$ such that for all $s \in \cS$, $f(s) \in X(s)
\otimes_{\sigma} H$. We equip $\cH_0$
with the inner product
$$\langle \delta_s \cdot \xi, \delta_t \cdot \eta \rangle = \delta_{s,t} \langle \xi, \eta \rangle  ,$$
for all $s,t \in \cS, \xi \in X(s) \otimes H, \eta \in
X(t) \otimes H$ (where the $\delta$'s on the left hand side are
Dirac deltas, the $\delta$ on the right hand side is Kronecker's
delta). Let $\cH$ be the completion of $\cH_0$ with
respect to this inner product. Note that
$$\cH \cong \oplus_{s \in \cS} X(s)\otimes H ,$$
but defining it as we did has a small notational advantage. We
define a family $\hat{T} = \{\hat{T}_s\}_{s \in \cS}$ of operators
on $\cH_0$ as follows. First, we define
$\hat{T}_0$ to be the identity. Now assume that $s>0$. If $t\in \cS$ and $t \ngeq s$, then we define $\hat{T}_s (\delta_t \cdot \xi ) = 0$ for all
$\xi \in X(t) \otimes_{\sigma} H$. And we define
\be\label{eq:def:hat} \hat{T}_s \left(\delta_t \cdot (x_{t-s}
\otimes x_s \otimes h) \right) = \delta_{t-s} \cdot
\left(x_{t-s}\otimes \tilde{T}_s (x_s \otimes h) \right)
\ee
if $t
\geq s > 0$. In \cite{Shalit07} we showed that $\hat{T} = \{\hat{T}_s\}_{s\in\cS}$ extends to a well defined semigroup of contractions on $\cH$.

Note that the adjoint of $\hat{T}$ is given by
$$\hat{T}_s \big(\delta_t \cdot x_t \otimes h \big) = \delta_{t+s}\cdot x_t \otimes \tilde{T}^*_s h ,$$
thus, if $T$ is a fully coisometric representation, then $\hat{T}$ is a semigroup of coisometries.

We summarize the construction in the following proposition.
\begin{proposition}\label{prop:technology}
Let $\cA$, $X$, $\cS$ and $(\sigma,T)$ be as above, and let
$$\cH = \oplus_{s \in \cS} X(s)\otimes_\sigma H.$$
There exists a contractive semigroup $\hat{T} = \{\hat{T}_s\}_{s\in\cS}$ on $\cH$ such that for all
$s\in\cS$, $x \in X(s)$ and $h\in H$,
$$\hat{T}_s \left(\delta_s \cdot x \otimes h \right) = T_s(x)h .$$
If $T$ is a fully coisometric representation, then $\hat{T}$ is a semigroup of coisometries.
\end{proposition}

\subsection{Isometric dilation of a fully coisometric representation}\label{subsec:coisoisodil}

For any $r = (r_1, \ldots, r_k) \in \mathbb{R}^k$, we denote $r_+
:= (\max\{r_1,0\}, \ldots, \max\{r_k,0\})$ and $r_- := r_+ - r$.
Throughout this section, $\cS$ will be a subsemigroup of $\Rp^k$
such that for all $s \in \cS - \cS$, both $s_+$ and $s_-$ are in
$\cS$. The semigroup that we are most interested in, namely
$\Rp^k$, satisfies this condition. For possible applications
discussed in Section \ref{subsec:k_tuple} we may need the
following theorem for $\mathbb{N}^k$, which also satisfies this
condition.
\begin{theorem}\label{thm:isoDilFC}
Let $\cS$ be as above, let $X = \{X(s)\}_{s \in \cS}$ be a product
system of unital $\cA$-correspondences over $\cS$, and let
$(\sigma,T)$ be a fully coisometric representation of $X$ on $H$,
with $\sigma$ unital. Then there exists a Hilbert space $K
\supseteq H$ and a minimal, fully coisometric and isometric
representation $(\rho,V)$ of $X$ on $K$, with $\rho$ unital, such
that
\begin{enumerate}
    \item $P_H$ commutes with $\rho (\cA)$, and $\rho(a) P_H = \sigma(a) P_H$, for all $a \in \cA$.
    \item\label{it:dilation} $P_H V_s(x)\big|_H = T_s(x)$ for all $s \in \cS$, $x \in X(s)$.
    \item\label{it:V*} $P_H V_s(x)\big|_{K \ominus H} = 0$ for all $s \in \cS$, $x \in X(s)$.
\end{enumerate}
If $\sigma$ is nondegenerate and $X$ is essential (that is, $\cA
X(s)$ is dense in $X(s)$ for all $s \in \cS$) then $\rho$ is also
nondegenerate. If $\cA$ is a $W^*$-algebra, $X$ is a product
system of $W^*$-correspondences and $(\sigma,T)$ is a
representation of $W^*$-correspondences, then $(\rho,V)$ is also a
representation of $W^*$-correspondences.
\end{theorem}
\begin{proof}
The proof is very similar to the proof of Proposition 3.2 in
\cite{Shalit07}, so we will not go into all the details whenever
they were taken care of in that paper. However, we note that there
are some essential differences between the situation at hand and
the one treated in \cite{Shalit07}.

Let $\cH = \oplus_{s \in \cS} X(s)\otimes_\sigma H$, and let
$\hat{T}$ be the semigroup of coisometries constructed in the
discussion preceding Proposition \ref{prop:technology}.%
%

Since $\hat{T}$ is a semigroup of coisometries, there exists a
minimal, \emph{regular} unitary dilation $W = \{W_s\}_{s\in \cS}$
of the semigroup $\{\hat{T}_s^*\}_{s\in \cS}$ on a Hilbert space
$\cK \supset \cH$ (this should be well known folklore, see
\cite{Shalit07b} for details). We denote $\hat{V}_s = W_s^*$. We
have for all $s \in \cS-\cS$
\be P_\cH \hat{V}_{s_+} \hat{V}_{s_-}^*
P_\cH = \hat{T}_{s_+} \hat{T}_{s_-}^* ,
\ee
Since the semigroup
$\hat{V}$ consists of commuting unitaries, and since commuting
unitaries doubly commute, we also have
\be\label{eq:reg_dil} P_\cH
\hat{V}_{s_-}^* \hat{V}_{s_+} P_\cH = \hat{T}_{s_+}
\hat{T}_{s_-}^* .
\ee
This triviality turns out to be crucial: it
will allow us to compute the inner products in $\cK$.

Introduce the Hilbert space $K$,
\bes K = \bigvee
\{\hat{V}_s(\delta_s \cdot(x \otimes h)) : s \in \cS, x \in X(s),
h \in H \} .
\ees
We consider $H$ as embedded in $K$ (or in $\cH$ or in $\cK$)
by the identification
\bes h \leftrightarrow \delta_0 \cdot(1
\otimes h) .
\ees
(This is where we use the fact that $\sigma$ is
unital). We turn to the definition of the representation $V$ of
$X$ on $K$. First, note that $\sigma(a)h$ is identified with
$\delta_0 \cdot 1 \otimes_{\sigma}\sigma(a)h = \delta_0 \cdot a
\otimes_{\sigma}h$. Next, we define a left action of $\cA$ on $\cH$ by
\bes
a \cdot (\delta_s \cdot x \otimes h) = \delta_s \cdot ax
\otimes h ,
\ees
for all $a \in \cA, s \in \cS, x \in X(s)$ and $h
\in H$. As we have explained in \cite{Shalit07b}, this
gives rise to a well defined a $*$-representation
that commutes
with $\hat{T}$: \bes a \hat{T}_s (\delta_t x_{t-s} \otimes x_s
\otimes h) = \delta_{t-s} a x_{t-s} \otimes T_s (x_s) h =
\hat{T}_s (\delta_t a x_{t-s} \otimes x_s \otimes h) . \ees
Taking adjoints shows that this left action commutes $\hat{T}_s^*$, ($s \in \cS$), as well.

We shall now define a representation $(\rho,V)$ of $X$ on $K$. First, we define
$\rho$ by the rule
\be \label{eq:V_e definition}
\rho(a) \hat{V}_s (\delta_s \cdot x_s \otimes h) = \hat{V}_s (\delta_s \cdot a x_s \otimes h) .
\ee
Using (\ref{eq:reg_dil}),
one shows that $\rho(a)$ extends to a bounded map on $K$. It then
follows by direct computation that $\rho$ is a $*$-representation.
Whe $(\sigma,T)$ is a representation of $W^*$-correspondences, we
also have to show that $\rho$ is a \emph{normal} representation.
Let $\{a_\gamma\} \subseteq {\rm ball}_1(\cA)$ be a net converging
in the weak operator topology to $a \in {\rm ball}_1(\cA)$. It is
known (for an outline of a proof, see \cite{MS03}) that the
mapping taking $b\in \cA$ to $b \otimes I_H \in B(X(s)
\otimes_\sigma H)$ is continuous in the ($\sigma$-)weak
topologies. Thus, for all $s \in \cS, x \in X(s)$ and $h \in H$,
$$a_\gamma  x \otimes h \longrightarrow a  x \otimes h $$
in the weak topology of $X(s) \otimes_\sigma H$. It follows that
$$\delta_s \cdot a_\gamma  x \otimes h \longrightarrow \delta_s \cdot a  x \otimes h $$
in the weak topology of $K$, so
$$\hat{V}_s \delta_s \cdot a_\gamma  x \otimes h \longrightarrow \hat{V}_s \delta_s \cdot a  x \otimes h $$
weakly. This implies that $\rho(a_\gamma) \rightarrow \rho(a)$ in the weak operator topology of $B(K)$, so $\rho$ is normal.

Note that $H$ reduces $\rho (A)$, and that $\rho(a) \big|_H =
\sigma(a) \big|_H$ (under the appropriate identifications).
Indeed, putting $t = 0$ in equation (\ref{eq:V_e definition})
gives \bes \rho(a) (\delta_0 \cdot 1 \otimes h) = \delta_0 \cdot a
\otimes h = \delta_0 \cdot 1  \otimes \sigma(a)h . \ees The
assertions regarding the unitality and nondegenracy of $\rho$ are
clear from the definitions.

We have completed the construction of $\rho$, and we proceed to
define the representation $V$ of $X$ on $K$. For $s
> 0$, we define $V_s$
by the rule \be \label{eq:definition V_s} V_s(x_s)
\hat{V}_t(\delta_t \cdot x_t \otimes h) = \hat{V}_{s+t}
(\delta_{s+t} \cdot x_s \otimes x_t \otimes h) . \ee
One has to use (\ref{eq:reg_dil}) to show that $V_s(x_s)$ can be extended to a well defined operator on $K$, but once that is done, it is easy to see that for all $s \in \cS$, $(\rho, V_s)$ is a covariant representation of $X(s)$ on
K. We now show that it is isometric. This computation is included so the reader has an opportunity to appreciate the role played by equation (\ref{eq:reg_dil}). Let $s,t,u \in \cS$, $x, y \in X(s)$, $x_t \in X(t)$, $x_u \in X(u)$ and $h,g \in H$. Then
\begin{align*}
\langle V_s(x)^* V_s(y) \hat{V}_t \delta_t \cdot x_t \otimes h, \hat{V}_u \delta_u \cdot x_u \otimes g \rangle
&= \langle  \hat{V}_{t+s} \delta_{t+s} \cdot y \otimes x_t \otimes h, \hat{V}_{u+s} \delta_{u+s} \cdot x \otimes x_u \otimes g \rangle \\
&= \langle  \hat{V}_{(t-u)_-}^* \hat{V}_{(t-u)_+} \delta_{t+s} \cdot y \otimes x_t \otimes h,  \delta_{u+s} \cdot x \otimes x_u \otimes g \rangle \\
(*)&= \langle  \hat{T}_{(t-u)_+} \hat{T}_{(t-u)_-}^* \delta_{t+s} \cdot y \otimes x_t \otimes h,  \delta_{u+s} \cdot x \otimes x_u \otimes g \rangle \\
&= \langle  \delta_{u+s} \cdot y \otimes \left(I \otimes \tilde{T}_{(t-u)_+}\right)\left(I \otimes \tilde{T}_{(t-u)_-}^* \right) \ldots \\
& \quad \ldots (x_t \otimes h),  \delta_{u+s} \cdot x \otimes x_u \otimes g \rangle \\
&= \langle  \delta_{u} \cdot \left(I \otimes \tilde{T}_{(t-u)_+}\right)\left(I \otimes \tilde{T}_{(t-u)_-}^* \right) \ldots \\
& \quad \ldots  (x_t \otimes h),  \delta_{u} \cdot \langle y, x\rangle x_u \otimes g \rangle \\
&= \langle  \hat{T}_{(t-u)_+} \hat{T}_{(t-u)_-}^* \delta_{t} \cdot x_t \otimes h ,  \delta_{u} \cdot \langle y, x\rangle x_u \otimes g \rangle \\
&= \langle  \hat{T}_{(t-u)_+} \hat{T}_{(t-u)_-}^* \delta_{t} \cdot \langle x, y\rangle x_t \otimes h ,  \delta_{u} \cdot \langle y, x\rangle x_u \otimes g \rangle \\
(*)&= \langle  \hat{V}_{(t-u)_-}^* \hat{V}_{(t-u)_+} \delta_{t} \cdot \langle x, y\rangle x_t \otimes h,  \delta_{u} \cdot  x_u \otimes g \rangle \\
&= \langle  \hat{V}_{t} \delta_{t} \cdot \langle x, y\rangle x_t \otimes h, \hat{V}_{u} \delta_{u} \cdot  x_u \otimes g \rangle \\
&= \langle \rho(\langle x,y \rangle) \hat{V}_t \delta_t \cdot x_t \otimes h, \hat{V}_u \delta_u \cdot x_u \otimes g \rangle .
\end{align*}
(The equations marked by (*) are where we use (\ref{eq:reg_dil}).)
This shows that $V_s(x)^* V_s(y) = \rho(\langle x,y \rangle)$, so
$(\rho,V)$ is indeed an isometric representation. To see that it
is fully coisometric, is is enough to show that for all $s\in\cS$,
$\tilde{V}_s$ is onto. It is clear that
$${\rm Im}(\tilde{V}_s) = \bigvee\{\hat{V}_{t+s} (\delta_{t+s} \cdot x_s \otimes x_t \otimes h) : t\in \cS, x_s \in X(s), x_t \in X(t), h\in H \} .$$
But if $t\in \cS$, $x_t \in X(t)$ and $h\in H$, then
\begin{align*}
\hat{V}_t (\delta_t \cdot x_t \otimes h)
&= \hat{V}_t \hat{V}_s \hat{V}_s^* (\delta_t \cdot x_t \otimes h) \\
(*)&= \hat{V}_t \hat{V}_s \hat{T}_s^* (\delta_t \cdot x_t \otimes h) \\
&= \hat{V}_{t+s}  (\delta_{t+s} \cdot x_t \otimes \tilde{T}_s^* h) \in {\rm Im}(\tilde{T}_s) ,
\end{align*}
where (*) is justified because $\hat{V}^*_s$ is an extension of
$\hat{T}^*_s$ (as is any unitary dilation of an isometry). This
shows that $\tilde{V}_s$ is onto, so it is a unitary, hence $V$ is
a fully coisometric.

We have yet to show that $V$ is a representation of product systems (that is, that the
semigroup property holds) and that it is in fact a dilation of $T$.

Let $h \in H$, $s,t,u \in \cS$, and let $x_s, x_t, x_u$ be in $X(s), X(t), X(u)$, respectively. Then
\begin{align*}
V_{s+t} (x_s \otimes x_t) \hat{V}_u (\delta_u \cdot x_u \otimes h)
& = \hat{V}_{s+t+u}(\delta_{s+t+u} \cdot x_s \otimes x_t \otimes x_u \otimes h) \\
& = V_s (x_s) \hat{V}_{t+u}(\delta_{t+u} \cdot x_t \otimes x_u \otimes h) \\
& = V_s (x_s) V_t (x_t) \hat{V}_{u}(\delta_{u} \cdot x_u \otimes h),
\end{align*}
so the semigroup property holds. Finally, let $s \in \cS, x \in X(s)$ and
$h = \delta_0 \cdot 1 \otimes h \in H$. We compute:
\begin{align*}
P_H V_s(x) \big|_H h &= P_H V_s(x) \delta_0 \cdot 1 \otimes h \\
&= P_H \hat{V}_s (\delta_s \cdot x \otimes h) \\
&= P_H P_{\cH} \hat{V}_s \big|_\cH (\delta_s \cdot x \otimes h) \\
&= P_H \hat{T}_s (\delta_s \cdot x \otimes h) \\
&= P_H (\delta_0 \cdot 1 \otimes T_s(x)h) = T_s(x)h .
\end{align*}
We remark that $V$ is already a minimal isometric dilation of $T$, because
\begin{align*}
K
&= \bigvee \{\hat{V}_s(\delta_s \cdot(x \otimes h)) : s \in \cS, x \in X(s),
h \in H \} \\
&= \bigvee \{V_s(x) (\delta_0 \cdot(1 \otimes h)) : s \in \cS, x \in X(s),
h \in H \} .
\end{align*}
Item \ref{it:V*} in the statement of the theorem follows as in Proposition 3.2, \cite{Shalit07}.
\end{proof}

\section{E$_0$-dilation of a strongly commuting pair of CP$_0$-maps}\label{sec:Edilation}

In this section we prove the main result of this paper: every pair
of strongly commuting CP$_0$-semigroups has a minimal
E$_0$-dilation. In the last two sections we worked out the two
main steps in the Muhly-Solel approach to dilation. In this
section we will put together these two steps and take care of the
remaining technicalities. It is convenient to begin by proving a
few technical lemmas. We then turn to prove the existence of the
dilation, and we close this section with a discussion of
minimality issues.

\subsection{Continuity of CP-semigroups}
\begin{lemma}\label{lem:semigroup}
Let $N$ be a von Neumann algebra, let $\cS$ be an abelian, cancellative semigroup with unit $0$, and let $X$ be a product system of $N$-correspondences over $\cS$. Let $W$ be completely contractive covariant representation of $X$ on a Hilbert space $G$, such that $W_0$ is unital. Then the family of maps
\bes
\Theta_s : a \mapsto \tilde{W}_s (I_{X(s)} \otimes a) \tilde{W}_s^* \,\, , \,\, a \in W_0 (N)',
\ees
is a semigroup of CP maps (indexed by $\cS$) on $W_0 (N)'$. Moreover, if $W$ is an isometric (a fully coisometric) representation, then $\Theta_s$ is a $*$-endomorphism (a unital map) for all $s\in\cS$.
\end{lemma}
\begin{proof}
By Proposition 2.21 in \cite{MS02}, $\{\Theta_s\}_{s\in\cS}$ is a family of contractive, normal, completely positive maps on $W_0(N)'$. Moreover, these maps are unital if $W$ is a fully coisometric representation, and they are $*$-endomorphisms if $W$ is an isometric representation. All that remains is to check that $\Theta = \{\Theta_s \}_{s\in\cS}$ satisfies the semigroup condition $\Theta_s \circ \Theta_t = \Theta_{s + t}$. Fix $a \in W_0 (N)'$. For all $s,t\in\cS$,
\begin{align*}
\Theta_s (\Theta_t (a))
&= \tilde{W}_s \left(I_{X(s)} \otimes \left(\tilde{W}_t (I_{X(t)} \otimes a) \tilde{W}_t^*\right)\right) \tilde{W}_s^* \\
&= \tilde{W}_s (I_{X(s)} \otimes \tilde{W}_t) (I_{X(s)}\otimes I_{X(t)} \otimes a)(I_{X(s)} \otimes \tilde{W}_t^*) \tilde{W}_s^* \\
&= \tilde{W}_{s + t} (U_{s,t}\otimes I_G)(I_{X(s)}\otimes I_{X(t)}\otimes a)(U_{s,t}^{-1}\otimes I_G)\tilde{W}_{s + t}^* \\
&= \tilde{W}_{s + t} (I_{X(s\cdot t)}\otimes a)\tilde{W}_{s + t}^* \\
&= \Theta_{s + t}(a) .
\end{align*}
Using the fact that $W_0$ is unital, we have
\begin{align*}
\Theta_0(a) h
&= \tilde{W_0} (I_N \otimes a) \tilde{W_0}^* h \\
&= \tilde{W_0} (I_N \otimes a) (I \otimes h) \\
&= W_0(I_N)ah \\
&= ah ,
\end{align*}
thus $\Theta_0(a) = a$ for all $a\in N$.
\end{proof}
\begin{lemma}\label{lem:continuity1}
Let $\{R_t\}_{t\geq0}$ and $\{S_t\}_{t\geq0}$ be two CP-semigroups
on $\cM \subseteq B(H)$, where $H$ is a separable Hilbert space.
Then the two parameter CP-semigroup $P$ defined by
$$P_{(s,t)} := R_s S_t$$
is a CP-semigroup, that is, for all $a \in \cM$, the map $\Rpt \ni
(s,t) \mapsto P_{(s,t)}(a)$ is weakly continuous. Moreover, $P$ is
\emph{jointly} continuous on $\Rpt \times \cM$, endowed with the
standard$\times$weak-operator topology.
\end{lemma}
\begin{proof}
Let $(s_n,t_n) \rightarrow (s,t)\in \Rpt$, and let $a_n
\rightarrow a \in \cM$. Then, by Proposition 4.1, \cite{MS02},
$S_{t_n}(a_n) \rightarrow S_t(a)$ in the weak operator topology.
By the same proposition used once more,
$$P_{(s_n,t_n)}(a_n) = R_{s_n}(S_{t_n}(a_n)) \rightarrow R_s(S_t(a)) = P_{(s,t)}(a)$$
where convergence is in the weak operator topology.
\end{proof}

The above lemma show that, given two CP$_0$-semigroups $\{R_t\}_{t\geq0}$ and
$\{S_t\}_{t\geq0}$, we can form a two-parameter CP$_0$-semigroup $\{P_{(s,t)}\}= \{R_s S_t\}_{s,t\geq0}$ which  satisfies
the natural continuity conditions. For the theorem below, we will need $P$ to satisfy a stronger type of continuity. This is the subject of the next two lemmas.

\begin{lemma}\label{lem:continuity2}
Let $\cS$ be a topological semigroup with unit $0$, and let $\{W_s\}_{s\in\cS}$ be a semigroup over $\cS$ of CP maps on a von Neumann algebra $\cR \subseteq B(H)$. Let $\cA\subseteq \cR$ be a sub $C^*$-algebra of $\cR$ such that for all $a \in \cA$,
$$W_s(a) \stackrel{WOT}{\longrightarrow} a$$
as $s \rightarrow 0$. Then  for all $a \in \cA$,
$$W_{t+s}(a) \stackrel{SOT}{\longrightarrow} W_t (a)$$
as $s \rightarrow 0$.
\end{lemma}
\begin{proof}
The proof is taken, almost word for word, from the proof of the first half of Proposition 4.1, \cite{MS02},
which addresses the case $\cS=\Rp$.

Let $a \in \cA$. It is enough to prove $W_{s}(a) \stackrel{SOT}{\longrightarrow} a$, as the result for $t$ other than $0$ follows from the normality of $W_t$ and from the semigroup property. Also, we may assume that $a$ is unitary. Let $h\in H$ be a unit vector. Then
$$\|W_s(a)h -ah \|^2 = \|W_s(a)h\|^2 -2 {\rm Re}\langle ah, W_s(a)h \rangle + \|a h\|^2.$$
To show that the right hand side converges to $0$ as $s \rightarrow 0$, it is enough to show that $\lim_{s\rightarrow 0}\|W_s(a)h\|^2 = \|a h\|^2 = 1$. But
\begin{align*}
1 &\geq \|W_s(a)h\| \\
&\geq |\langle W_s (a)h, ah \rangle| \stackrel{s \rightarrow 0}{\longrightarrow} |\langle a h, ah \rangle | = 1.
\end{align*}
This completes the proof.
\end{proof}

\begin{lemma}\label{lem:continuity3}
Let $\Theta = \{\Theta_t\}_{t\geq 0}$ be a CP-semigroup on $\cM \subseteq B(H)$, where $H$ is a separable Hilbert space. Then $\Theta$ is jointly strongly continuous, that is, for all $h\in H$,
the map
$$(t,a) \mapsto \Theta_t(a)h$$
is continuous in the standard$\times$strong-operator topology.
\end{lemma}
\begin{proof}
First, assume that $\Theta$ is an E-semigroup. Let $(t_n,a_n)
\rightarrow (t,a)$ in the standard$\times$strong-operator topology
in $\Rp \times \cM$, and $h\in H$. \bes \|\Theta_{t_n} (a_n) h -
\Theta_t(a)h \|^2 = \|\Theta_{t_n} (a_n) h\|^2 - 2 {\rm Re}\langle
\Theta_{t_n} (a_n) h,\Theta_{t} (a) h \rangle + \|\Theta_{t} (a)
h\|^2 , \ees since $\Theta$ is continuous in the
standard$\times$weak-operator topology, it is enough to show that
$\|\Theta_{t_n} (a_n) h\|^2 \rightarrow \|\Theta_{t} (a) h\|^2$.
But \bes \| \Theta_{t_n} (a_n) h\|^2 = \langle \Theta_{t_n} (a_n^*
a_n) h, h \rangle \rightarrow \langle \Theta_{t} (a^* a) h, h
\rangle = \| \Theta_{t} (a) h\|^2 , \ees because $a_n^* a_n
\rightarrow a^* a$ in the weak-operator topology, and $\Theta$ is
jointly continuous with respect to this topology. Thus $\Theta$ is
also jointly continuous with respect to the strong-operator
topology.

Now let $\Theta$ be an arbitrary CP-semigroup, and let $(K,u,\cR,\alpha)$ be an E-dilation of $\Theta$. Then
for all $a \in \cM, t \in \Rp$,
$$\Theta_t(a) = u^* \alpha_t(u a u^* ) u ,$$
whence $\Theta$ inherits the required type of joint continuity from $\alpha$.
\end{proof}

From the above lemma one immediately obtains:

\begin{proposition}\label{prop:continuity4}
Let $\{R_t\}_{t\geq0}$ and
$\{S_t\}_{t\geq0}$ be two CP-semigroups on $\cM \subseteq B(H)$, where $H$ is a separable Hilbert space. Then the
two parameter CP-semigroup $P$ defined by
$$P_{(s,t)} := R_s S_t$$
is strongly continuous, that is, for all $a \in \cM$, the map $\Rpt \ni (s,t) \mapsto P_{(s,t)}(a)$ is strongly continuous. Moreover, $P$ is \emph{jointly} continuous on $\Rpt \times \cM$, endowed with the standard$\times$strong-operator topology.

\end{proposition}
%
%

\subsection{The existence of an E$_0$-dilation}
We have now gathered enough tools to prove our main
result.
\begin{theorem}\label{thm:scudil}
Let $\{R_t\}_{t\geq0}$ and $\{S_t\}_{t\geq0}$ be two strongly
commuting CP$_0$-semigroups on a von Neumann algebra $\cM\subseteq
B(H)$, where $H$ is a separable Hilbert space. Then the two
parameter CP$_0$-semigroup $P$ defined by
$$P_{(s,t)} := R_s S_t$$
has a minimal E$_0$-dilation $(K,u,\cR,\alpha)$. Moreover, $K$ is separable.
\end{theorem}
\begin{proof}
We split the proof into the following steps:
\begin{enumerate}
\item Existence of a $*$-endomorphic dilation $(K,u,\cR,\alpha)$ for $(\cM,P)$.
\item Minimality of the dilation.
\item Continuity of $\alpha$ on $\cM$.
\item Separablity of $K$.
\item Continuity of $\alpha$.
\end{enumerate}

{\bf Step 1: Existence of a $*$-endomorphic dilation}

Let $X$ and $T$ be the product system (of $\cM'$-correspondences) and the fully coisometric
product system representation given by Theorem \ref{thm:rep_SC}. By Theorem \ref{thm:isoDilFC}, there is
a covariant isometric and fully coisometric representation $(\rho,V)$ of $X$ on some Hilbert space $K \supseteq H$, with $\rho$ unital.
Put $\tilde{\cR} = \rho(M')'$, and let $u$ be the isometric inclusion $H \rightarrow K$. Note that, since $u H$ reduces $\rho$, $p:= u u^* \in \tilde{\cR}$. We define a semigroup
$\tilde{\alpha}=\{\tilde{\alpha}_s\}_{s\in\Rpt}$ by
\bes
\tilde{\alpha}_s(b) = \tilde{V}_s(I \otimes b)\tilde{V}_s^* \,\, , \,\, s\in\Rpt ,
b\in \tilde{\cR}.
\ees
By Lemma \ref{lem:semigroup} above, $\tilde{\alpha}$ is a semigroup of unital, normal
$*$-endomorphisms of $\tilde{\cR}$. The (first part of the) proof of Theorem
2.24 in \cite{MS02} works in this situation as well, and shows
that
\be\label{eq:M=pRp}
\cM = u^* \tilde{\cR} u
\ee
and that
\be\label{eq:dilationeq}
P_s (u^* b u) = u^* \tilde{\alpha}_s(b) u
\ee
for all
$b\in \tilde{\cR}$, $s \in \Rpt$. Note that we \emph{cannot} use that theorem directly, because
for fixed $s\in\cS$, $X(s)$ is not necessarily the identity representation of $P_s$. For the sake of completeness, we repeat the argument (with some changes).

By Theorem \ref{thm:isoDilFC},
for all $a \in \cM'$, $u^* \rho(a) u = \sigma (a)$, and by definition, $\sigma(a) = a$, thus
\bes\label{eq:corner}
u^* \tilde{\cR} u = u^*\rho(\cM')' u = (u^*\rho(\cM')u)' = (\cM')' = \cM,
\ees
where the $\subseteq$ part of the second equality follows from the fact that $uH$ reduces $\rho(\cM')$. This establishes (\ref{eq:M=pRp}), which allows us to make the identification $\cM = p \tilde{\cR} p \subseteq \tilde{\cR}$. To obtain (\ref{eq:dilationeq}), we fix $s \in \Rpt$ and $b \in \tilde{\cR}$, and we compute
\begin{align*}
P_s(u^* b u)
&= \tT_s(I \otimes u^* bu) \tT_s^*  \\
(*)&= u^*\tV_s(I \otimes u)(I \otimes u^* bu)(I \otimes u^*) \tV_s^* u  \\
(**)&= u^*\tV_s(I \otimes  b) \tV_s^* u  \\
&= u^* \tilde{\alpha}_s (b) u .
\end{align*}
The equalities marked by (*) and (**) are justified by items \ref{it:dilation} and \ref{it:V*} of Theorem \ref{thm:isoDilFC}, respectively. Equation (\ref{eq:dilationeq}) implies that $p$ is a coinvariant
projection. Since $\tilde{\alpha}$ is unital, we have $\tilde{\alpha}_t(p) \geq p$ for all $t \in \Rpt$, that is, $p$  is an increasing projection.

Even though we started out with a minimal isometric representation
$V$ of $T$, we cannot show that $\tilde{\alpha}$ is a
minimal dilation of $P$. We define
\be\label{eq:redefine R}
\cR =
W^* \left(\bigcup_{t\in\Rpt}\tilde{\alpha}_t(\cM) \right) .
\ee
This von Neumann algebra is invariant under $\tilde{\alpha}$, and we denote $\alpha = \tilde{\alpha}|_{\cR}$.
Now it is immediate that $(p,\cR,\alpha)$ is a $*$-endomorphic dilation of $(\cM,P)$. Indeed, for all $b\in\cR$ and all $t \in\Rpt$,
$$p \alpha_t (b) p = p \tilde{\alpha}_t(b) p = P_t (pbp) ,$$
because $(p,\tilde{\cR},\tilde{\alpha})$ is a dilation of $(\cM,P)$. It is also clear that $\cM = p\cR p$.

The only issue left to handle is the continuity of $\alpha$. We now define two one-parameter semigroups on $\cR$: $\beta =
\{\beta_t\}_{t\geq 0}$ and $\gamma = \{\gamma_t\}_{t\geq 0}$ by $\beta_t = {\alpha}_{(t,0)}$ and $\gamma_t =
{\alpha}_{(0,t)}$. Clearly, $\beta$ and $\gamma$ are semigroups of normal, unital $*$-endomorphisms of $\cR$. If we show that $K$ is separable, then by Lemma \ref{lem:continuity1}, once we show that $\beta$ and $\gamma$ are E$_0$-semigroups -- that is, possess the required weak continuity -- then we have shown that $\alpha$ is an E$_0$-semigroup. The rest of the proof is dedicated to showing that $\beta$ and $\gamma$ are E$_0$-semigroups and that $K$ is separable. But before we do that, we must show that the dilation is minimal, and, in fact, a bit more.

{\bf Step 2: Minimality of the dilation}

What we really need to prove is that
\be\label{eq:partition1}
K = \bigvee \alpha_{(s_m,t_n)}(\cM) \alpha_{(s_m,t_{n-1})}(\cM) \cdots \alpha_{(s_m,t_1)}(\cM) \alpha_{(s_{m},0)}(\cM) \alpha_{(s_{m-1},0)}(\cM) \cdots \alpha_{(s_1,0)}(\cM) H
\ee
where in the right hand side of the above expression we run over all strictly positive pairs $(s,t)\in\Rpt$ and all partitions $\{0=s_0 < \ldots < s_m=s\}$ and $\{0=t_0 < \ldots < t_n = t \}$ of $[0,s]$ and $[0,t]$. We shall also need an analog of (\ref{eq:partition1}) with the roles of the first and second ``time variables" of $\alpha$ replaced, but since the proof is very similar we shall not prove it.

Recall that
$$K = \bigvee \left\{V_{(s,t)}(X(s,t))H : (s,t)\in\Rpt \right\} .$$
Thus, it suffices to show that for a fixed $(s,t)\in\Rpt$,
\be\label{eq:partition2}
V_{(s,t)}(X(s,t))H = \bigvee \alpha_{(s_m,t_n)}(\cM) \cdots \alpha_{(s_m,t_1)}(\cM) \alpha_{(s_{m},0)}(\cM) \alpha_{(s_{m-1},0)}(\cM) \cdots \alpha_{(s_1,0)}(\cM) H
\ee
where in the right hand side of the above expression we run over all  partitions $\{0=s_0 < \ldots < s_m=s\}$ and $\{0=t_0 < \ldots < t_n = t \}$ of $[0,s]$ and $[0,t]$.

To show that we can consider only $s$ and $t$ strictly positive, we note that if $u,v \in\Rpt$, then
\begin{align*}
V_u(X(u))H
&= \tilde{V}_u(I_{X(u)}\otimes \tilde{V}_v)(I_{X(u)}\otimes \tilde{V}_v^*)(X(u) \otimes H) \\
&= \tilde{V}_u(I_{X(u)}\otimes \tilde{V}_v)(I_{X(u)}\otimes \tilde{T}_v^*)(X(u) \otimes H) \\
&= \tilde{V}_{u+v}(X(u)\otimes \tilde{T}_v^* H) \\
&\subseteq V_{u+v}(X(u+v))H .
\end{align*}

We now turn to establish (\ref{eq:partition2}). Recall the notation and constructions of Subsections \ref{subsec:des_MS} and \ref{subsec:repvia}.
$$X(s,t) := E(s) \otimes F(t) ,$$
and
$$T_{(s,t)}(\xi \otimes \eta) := T_s^E (\xi) T_t^F (\eta) ,$$
where $(E,T^E)$ and $(F,T^F)$ are the product systems and representations  representing $R$ and $S$ via Muhly and Solel's construction as described in \ref{subsec:des_MS}. By Lemma 4.3 (2) of \cite{MS02}, for all $r > 0$,
$$\bigvee \{(I_{E(r)} \otimes a ) (\tilde{T}_r^E)^* h : a \in \cM, h \in H \} = \cE_r \otimes_{\cM'} H,$$
where $\cE_r = \cL_M (H,H^R_\mfp)$ with the partition $\mfp = \{0=r_0 < r_1 = r\}$. Similarly,
$$\bigvee \{(I_{F(r)} \otimes a ) (\tilde{T}_r^F)^* h : a \in \cM, h \in H \} = \cF_r \otimes_{\cM'} H .$$
Fix $s,t>0$. Under the obvious identifications,  if we go over all the partitions $\{0=s_0 < \ldots < s_m=s\}$ and $\{0=t_0 < \ldots < t_n = t \}$ of $[0,s]$ and $[0,t]$, the collection of correspondences
$$\cE_{s_1} \otimes \cE_{s_2-s_1} \otimes \cdots \otimes \cE_{s_m-s_{m-1}} \otimes \cF_{t_1} \otimes \cdots \otimes \cF_{t_n - t_{n-1}}$$
is dense in $X(s,t)$. Using Lemma 4.3 (2) of \cite{MS02} repeatedly, we obtain
\begin{align*}
& \alpha_{(s_m,t_n)}(\cM) \cdots \alpha_{(s_m,t_1)}(\cM) \alpha_{(s_{m},0)}(\cM) \alpha_{(s_{m-1},0)}(\cM) \cdots \alpha_{(s_1,0)}(\cM) H \\
&= \alpha_{(s_m,t_n)}(\cM) \cdots \tilde{V}_{(s_1,0)}(I_{(s_1,0)} \otimes \cM)\tilde{V}_{(s_1,0)}^* H \\
&= \alpha_{(s_m,t_n)}(\cM) \cdots \tilde{V}_{(s_1,0)} (I_{(s_1,0)} \otimes \cM)(\tilde{T}^E_{s_1})^* H \\
&= \alpha_{(s_m,t_n)}(\cM) \cdots \tilde{V}_{(s_2,0)}(I_{(s_2,0)} \otimes \cM)\tilde{V}_{(s_2,0)}^* \tilde{V}_{(s_1,0)} (\cE_{s_1} \otimes H).
\end{align*}
But
$$\tilde{V}_{(s_2,0)}^* \tilde{V}_{(s_1,0)} = (I_{(s_1,0)}\otimes \tilde{V}_{(s_2-s_1,0)}^*)\tilde{V}_{(s_1,0)}^*\tilde{V}_{(s_1,0)} = (I_{(s_1,0)}\otimes \tilde{V}_{(s_2-s_1,0)}^*) ,$$
so we get
\begin{align*}
& \alpha_{(s_m,t_n)}(\cM) \cdots \alpha_{(s_m,t_1)}(\cM) \alpha_{(s_{m},0)}(\cM) \alpha_{(s_{m-1},0)}(\cM) \cdots \alpha_{(s_1,0)}(\cM) H \\
&= \alpha_{(s_m,t_n)}(\cM) \cdots \tilde{V}_{(s_2,0)}(I_{(s_2,0)} \otimes \cM)(I_{(s_1,0)}\otimes \tilde{V}_{(s_2-s_1,0)}^*) (\cE_{s_1} \otimes H) \\
&= \alpha_{(s_m,t_n)}(\cM) \cdots \tilde{V}_{(s_2,0)}(\cE_{s_1} \otimes \cE_{s_2-s_1} \otimes H).
\end{align*}
Continuing this way, we see that
\begin{align*}
& \alpha_{(s_m,t_n)}(\cM) \cdots \alpha_{(s_m,t_1)}(\cM) \alpha_{(s_{m},0)}(\cM) \alpha_{(s_{m-1},0)}(\cM) \cdots \alpha_{(s_1,0)}(\cM) H \\
&= {V}_{(s,t)} (\cE_{s_1} \otimes \cE_{s_2-s_1} \otimes \cdots \otimes \cE_{s_m-s_{m-1}} \otimes \cF_{t_1} \otimes \cdots \otimes \cF_{t_n - t_{n-1}})H .
\end{align*}
Since this computation works for any partition of $[0,s]$ and $[0,t]$, we have (\ref{eq:partition2}). This, in turn, implies (\ref{eq:partition1}), which is what we have been after.

Now it is a simple matter to show that $(p,\cR,\alpha)$ is a minimal dilation of $(\cM,P)$. First, note that by (\ref{eq:partition1})
$$K = \left[\cR p K \right].$$
In light of (\ref{eq:redefine R}), Definitions \ref{def:min_dil} and \ref{def:min_dil_Arv} and Proposition
\ref{prop:equiv_def_min}, we have to show that the central support of $p$ in $\cR$ is $I_K$. But this follows by a standard (and short) argument, which we omit.

{\bf Step 3: Continuity of $\beta$ and $\gamma$ on $\cM$}

We shall now show that function $\Rp \ni t \mapsto \beta_t(a)$ is strongly continuous from the right for each $a \in \cA := C^* \left(\bigcup_{t\in\Rpt} {\alpha}_t(M)\right)$. Of course, the same is true for $\gamma$ as well.

Let $k_1 = \sum_i {\alpha}_{s_i}(m_i)h_i$ and $k_2 = \sum_j
{\alpha}_{t_j}(n_j)g_j$  be in $K$, where $s_i = (s_i^1, s_i^2) , t_j = (t_j^1, t_j^2) \in \Rpt$,
$m_i,n_j \in \cR$ and $h_i,g_j \in H$. By (\ref{eq:partition1}), we may consider only $s_i^1, t_j^1 >0$. Take $a \in \cM$ and $t > 0$. For the following computations, we may assume that $k_1$ and $k_2$ are
given by finite sums, and we take $t <
\min\{t_j^1,s_i^1\}_{i,j}$. We will abuse notation a bit by denoting $(t,0)$ by $t$. Now compute:
\begin{align*}
\langle{\beta}_t(a)k_1,k_2\rangle
&= \sum_{i,j}\langle{\alpha}_t(a){\alpha}_{s_i}(m_i)h_i,{\alpha}_{t_j}(n_j)g_j\rangle \\
&= \sum_{i,j}\langle {\alpha}_{t_j}(n_j^*) {\alpha}_{t}(a){\alpha}_{s_i}(m_i)h_i,g_j\rangle \\
&= \sum_{i,j}\langle {\alpha}_{t}\left({\alpha}_{t_j-t}(n_j^*) a {\alpha}_{s_i-t}(m_i)\right)h_i,g_j\rangle \\
&= \sum_{i,j}\langle P_{t}\left(p{\alpha}_{t_j-t}(n_j^*) pap{\alpha}_{s_i-t}(m_i)p\right)h_i,g_j\rangle \\
&= \sum_{i,j}\langle P_{t}\left(P_{t_j-t}(p n_j^* p) a P_{s_i-t}(p m_i p)\right)h_i,g_j\rangle \\
&\stackrel{t\rightarrow 0}{\longrightarrow} \sum_{i,j}\langle P_{t_j}(p n_j^* p) a P_{s_i}(p m_i p) h_i,g_j\rangle \\
&= \sum_{i,j}\langle a {\alpha}_{s_i}(m_i)h_i,{\alpha}_{t_j}(n_j) g_j\rangle \\
&= \langle a k_1,k_2\rangle ,
\end{align*}
where we have made use of the joint strong continuity of $P$ (Proposition \ref{prop:continuity4}). This implies that for all $a\in \cM$,  ${\alpha}_t(a) \rightarrow a$ weakly as $t
\rightarrow 0$.
It follows from Lemma \ref{lem:continuity2} that $\beta$ is strongly right continuous on $\bigcup_{t\in\Rpt} {\alpha}_t(\cM)$, whence it is also strongly right continuous on $\cA := C^* \left(\bigcup_{t\in\Rpt} {\alpha}_t(M)\right)$.

{\bf Step 4: Separability of $K$}

As we have already noted in Step 2, from (\ref{eq:partition1}) it follows that
$$K = \bigvee\{\alpha_{u_1}(a_1) \cdots \alpha_{u_k}(a_k) h: u_i \in \Rpt, a_i \in \cM, h \in H \}.$$
We define
$$K_0 = \bigvee\{\gamma_{t_1}(\beta_{s_1}((a_1)) \cdots \gamma_{t_k}(\beta_{s_k}((a_k)) h: s_i,t_i \in \mathbb{Q}_+, a_i \in \cM, h \in H \},$$
and
$$K_1 = \bigvee\{\gamma_{t_1}(\beta_{s_1}((a_1)) \cdots \gamma_{t_k}(\beta_{s_k}((a_k)) h: s_i \in \Rp,t_i \in \mathbb{Q}_+, a_i \in \cM, h \in H \}.$$
$K_0$ is clearly separable. Because of the normality of $\gamma$, the strong right continuity of $\beta$ on $\cM$ and the fact that multiplication is strongly continuous on bounded subsets of $\cR$, we can assert that $K_0 = K_1$, thus $K_1$ is separable. Now from the strong right continuity of $\gamma$ on $\cA$ and the continuity of multiplication, we see that $K = K_1$, whence it is separable.

{\bf Step 5: Continuity of $\alpha$}

Recall that all that we have left to show is that $\beta$ and $\gamma$ possess the desired weak continuity. We shall concentrate on $\beta$.

A short summary of the situation: we have a semigroup $\beta$ of normal, unital $*$-endomorphisms defined on a von Neumann algebra $\cR$ (which acts on a \emph{separable} Hilbert space $K$), and there is a weakly dense $C^*$-algebra $\cA \subseteq \cR$ such that for all $a\in \cA, k\in K$, the function $\Rp \ni \tau \mapsto \beta_\tau (a) k  \in K$ is right continuous.
From this, we want to conclude that for all $b \in \cR$, and all $k_1, k_2 \in K$, the map
$$ \tau \mapsto \langle \beta_\tau (b) k_1, k_2 \rangle$$
is continuous. This problem was already handled by Arveson in
\cite{Arv03} and by Muhly and Solel in \cite{MS02}. For
completeness, we give some shortened variant of their arguments.

For every $b \in \cR$, there is a sequence $\{a_n\}$ in $\cA$
weakly converging to $b$. Thus, for every $b\in \cR$ and every
$k_1, k_2, \in K$, the function $ \tau \mapsto \langle \beta_\tau
(b) k_1, k_2 \rangle$ is the pointwise limit of the sequence of
right continuous functions $ \tau \mapsto \langle \beta_\tau (a_n)
k_1, k_2 \rangle$,  so it is measurable. It now follows from
Proposition 2.3.1 in \cite{Arv03} (which, in turn, follows from
well known results in the theory of operator semigroups) that
$\beta$ is an E$_0$-semigroup.
\end{proof}

By Proposition \ref{prop:SC_H_finite}, if $H$ is a finite dimensional Hilbert space, then every pair of commuting CP-semigroups on $B(H)$ commutes strongly. Denote by $M_n(\mathbb{C})$ the algebra of $n\times n$ complex matrices. We have the following corollary.
\begin{corollary}\label{cor:M_n(C)}
Every two parameter CP$_0$-semigroup on $M_n(\mathbb{C})$ has an E$_0$-dilation.
\end{corollary}

Loosely speaking, the whole point of dilation theory is to present a certain object as part of a simpler, better understood object. Theorem \ref{thm:scudil} tells us that a two-parameter CP$_0$-semigroup can always be dilated to a two parameter E$_0$-semigroup. Certainly, E$_0$-semigroups are a very special case of CP$_0$-semigroups, so we have indeed made the situation simpler. But did we really? Perhaps $P$ (the CP$_0$-semigroup) was acting on a very simple kind of von Neumann algebra, but now $\alpha$ (the dilation) is acting on a very complicated one? Actually, we did not say much about the structure of $\cR$ (the dilating algebra). In this context, we have the following partial, but quite satisfying, result.
\begin{proposition}
If $\cM = B(H)$, then $\cR = B(K)$.
\end{proposition}
\begin{proof}
Let $q\in B(K)$ be a projection in $\cR'$. In particular, $pq = qp = pqp$, so $qp$ is a projection
$B(H)$ which commutes with $B(H)$, thus $qp$ is either $0$ or $I_H$.

If it is $0$ then for all $t_i\in\Rpt, m_i\in \cM, h \in H$,
$$q \alpha_{t_1}(m_1)\cdots\alpha_{t_k}(m_k)h = \alpha_{t_1}(m_1)\cdots\alpha_{t_k}(m_k) qp h = 0 ,$$
so $qK = 0$ and $q=0$.

If $qp = I_H$ then for all $0<t_i\in\Rpt, m_i\in \cM, h \in H$,
\begin{align*}
q \alpha_{t_1}(m_1)\cdots\alpha_{t_k}(m_k)h
&= \alpha_{t_1}(m_1)\cdots\alpha_{t_k}(m_k) qp h \\
&= \alpha_{t_1}(m_1)\cdots\alpha_{t_k}(m_k)h ,
\end{align*}
so $qK = K$ and $q=I_K$. We see that the only projections in $\cR'$ are $0$ and $I_K$, so
$\cR' = \mathbb{C}\cdot I_K$, hence $\cR = \cR'' = B(K)$.
\end{proof}

\section{Prospects for further results}\label{sec:prospects}
In the previous section we proved the main result of this paper,
Theorem \ref{thm:scudil}, which says that every pair of strongly
commuting CP$_0$-semigroups has an E$_0$-dilation. In fact, the
only place where strong commutativity was used was in showing that
the CP$_0$-semigroup at hand could be represented by a product
system representation as in the following equation
\be\label{eq:rep_rep1}
\Theta_s (a) = \tilde{T_s} \left(I_{X(s)}
\otimes a \right) \tilde{T_s}^* .
\ee
Furthermore, in light of our
dilation result from Subsection \ref{subsec:coisoisodil}, Theorem
\ref{thm:isoDilFC}, we see that given a subsemigroup $\cS\subseteq
\mathbb{R}^k$ such that for all $s \in \cS$, $s_-,s_+ \in \cS$,
and a CP$_0$-semigroup $\Theta = \{\Theta_s\}_{s\in\cS}$ acting on
a von Neumann algebra $\cM \subseteq B(H)$, ($H$ separable), an
E$_0$-dilation of $\Theta$ can be constructed if we are able to
find a product system of $\cM'$-correspondences $X$ over $\cS$ and
a fully coisometric product system representation $T$ of $X$ on
$H$ fulfilling (\ref{eq:rep_rep1}). In this section we use this
observation to dilate a CP$_0$-semigroup over $\mathbb{N}\times
\Rp$ which does not satisfy strong commutation.

\subsection{Example: E$_0$-dilation of a CP$_0$-semigroup over
$\mathbb{N}\times \Rp$ - without strong commutation}\label{subsec:NSC}
\label{subsec:NtimesR} Let $H = \mathbb{C} \oplus L^2(0,\infty)$.
Denote by $U$ the left-shift semigroup on $L^2(0,\infty)$ given by
$$(U_t f)(s) = f(t+s).$$
Let $S_t = 1 \oplus U_t$, and define a CP$_0$-semigroup $\Phi$ on
$B(H)$ by
$$\Phi_t(a) = S_t a S_t^*.$$
Next, define $k=1\oplus 0 \in H$, and define the CP map $\Theta$
by
$$\Theta (a) = \langle a k,k \rangle I_H \,\, , \,\, a \in B(H) .$$
Peeking into Example 5.5 in \cite{S06} one sees that for all $t\in
\Rp$, $\Theta$ and $\Phi_t$ commute but not strongly. However, we
shall show that the CP$_0$-semigroup $\Psi =
\{\Psi_{n,t}\}_{(n,t)\in\mathbb{N}\times \Rp}$  defined by
$\Psi_{n,t} = \Theta^n \circ \Phi_t$ has an E$_0$-dilation. In
light of the opening remarks of this section, all we have to do is
construct an appropriate product system representation.

Let $\{e_i\}_{i=1}^\infty$ be an orthonormal basis for
$L^2(0,\infty)$, and set $e_0=k$. Define $E_{i,0}$ to be the
infinite square matrix indexed by $I = \{0,1,2,\ldots\}$ having
$1$ in the $i$th rowú $0$th column, and zeros elsewhere. Abusing notation
slightly, we let $E_{i,0}$ denote also the operator that this
matrix represents with respect to the basis $\cE =
\{e_0\}_{i=0}^\infty$, namely, the rank one operator $e_i \otimes
e_0^*$. We note that
$$\Theta(a) = \sum_{i\in I} E_{i,0}aE_{i,0}^* .$$

If $(a_{i,j}(t))_{i,j\in I}$ is the matrix representing $S_t$ with
respect to $\cE$, then we have
$$S_t e_j = \sum_{i\in I}a_{i,j}(t) e_i,$$
thus
$$S_t E_{j,0} = \sum_{i\in I}a_{i,j}(t) E_{i,0} = \sum_{i\in I}a_{i,j}(t) E_{i,0} S_t.$$
The matrix function $a(t)$ is a semigroup of coisometric matrices,
so there is a semigroup of unitary matrices $\{u(t)\}_{t\geq 0}$
indexed by $I \cup I'$, where $I'$ is another copy of $I$, such
that the $I$-$I$ block in $u(t)$ is equal to $a(t)$, and the
$I$-$I'$ block in $u(t)$ is $0$ ($u$ is simply the matrix
representation of the minimal isometric dilation of the semigroup
$S$, which is unitary, because $a(t)$ is coisometric). We now define a family $\{T_i\}_{i\in I\cup I'}$ of
operators on $H$ by $T_i = E_{i,0}$ when $i\in I$ and $T_i = 0$
when $i\in I'$. Because of the block structure that $u(t)$
possesses, we have for all $t\geq 0$ \be\label{unicomm}S_t T_j =
\sum_{i \in I\cup I'} u_{i,j}(t)T_i S_t.\ee We shall now construct
a product system of Hilbert spaces over $\mathbb{N}\times \Rp$.
Let $E = \ell^2(I\cup I')$, and put $E(n) = E^{\otimes n}$. We fix
an orthonormal basis $\cF=\{f_{i}\}_{i\in I\cup I'}$ in $E$. Also,
let $F$ be the trivial product system, that is, the product system with $F(t) = \mathbb{C}$ for all $t\in\Rp$ and
the obvious multiplication. For all
$n\in \mathbb{N}$ and all $t\in\Rp$, we define
$$X(n,t) = E(n) \otimes F(t).$$ To make $X =
\{X(n,t)\}_{(n,t)\in\mathbb{N}\times \Rp}$ into a product system,
we must define unitaries
$$U_{(m,s)(n,t)} : X(m,s)\otimes X(n,t) \rightarrow X(m+n,s+t)$$
that are associative in the sense of equation
(\ref{eq:assoc_prod}). This is where $u$ comes in. If $\lambda \in
F(s), \mu \in F(t)$, we define
$$U_{(1,s)(1,t)} (f_i \otimes \lambda)\otimes (f_j \otimes \mu) =
\sum_{k\in I\cup I'}u_{k,j}(t)\, f_i \otimes f_k \otimes \lambda \mu
,$$ and we continue this map to all of $X$. Let $k,m,n \in
\mathbb{N}$, and $s,t,u \in \Rp$. We have to show that
$$U_{(k,s)(m+n,t+u)}(I \otimes U_{(m,t)(n,u)}) = U_{(k+m,s+t)(n,u)}(U_{(k,s)(m,t)}\otimes I).$$
We shall operate with both sides on a typical element of the form
$$f_{i_1}\otimes \cdots \otimes f_{i_k}\otimes \lambda \otimes f_{j_1}\otimes
\cdots f_{j_m}\otimes \mu \otimes f_{l_1}\otimes \cdots \otimes
f_{l_n}\otimes \nu,$$ where $\lambda \in F(s)$, $\mu \in F(t)$ and
$\nu \in F(u)$. Operating first with $(I \otimes U_{(m,t)(n,u)})$
we get
$$\sum_{l_1',\ldots, l_n'} u_{l_1',l_1}(t) \cdots u_{l_n',l_n}(t) f_{i_1}\otimes \cdots \otimes f_{i_k}\otimes \lambda \otimes f_{j_1}\otimes
\cdots f_{j_m}\otimes f_{l_1'}\otimes \cdots \otimes
f_{l_n'}\otimes \mu \nu ,$$ and following with an application of
$U_{(k,s)(m+n,t+u)}$ we get
\begin{align*} \sum_{l_1',\ldots,
l_n'} u_{l_1',l_1}(t) \cdots u_{l_n',l_n}(t) \sum_{j_1',\ldots,
j_m'} u_{j_1',j_1}(s) \cdots u_{j_m',j_m}(s)
\sum_{l_1'',\ldots,l_n''} u_{l_1'',l_1'}(s) \cdots
u_{l_n'',l_n'}(s)\\
 f_{i_1}\otimes \cdots \otimes f_{i_k}
\otimes f_{j_1'}\otimes \cdots f_{j_m'}\otimes f_{l_1''}\otimes
\cdots \otimes f_{l_n''}\otimes \lambda \mu \nu
\end{align*}
which is
\begin{align*}  \sum_{j_1',\ldots,
j_m'} u_{j_1',j_1}(s) \cdots u_{j_m',j_m}(s)
\sum_{l_1'',\ldots,l_n''} u_{l_1'',l_1}(s+t) \cdots
u_{l_n'',l_n}(s+t)\\
 f_{i_1}\otimes \cdots \otimes f_{i_k}
\otimes f_{j_1'}\otimes \cdots f_{j_m'}\otimes f_{l_1''}\otimes
\cdots \otimes f_{l_n''}\otimes \lambda \mu \nu
\end{align*}
because $u$ is a semigroup. On the other hand applying first
$(U_{(k,s)(m,t)}\otimes I)$ we get
$$\sum_{j_1',\ldots,
j_m'} u_{j_1',j_1}(s) \cdots u_{j_m',j_m}(s)  f_{i_1}\otimes
\cdots \otimes f_{i_k} \otimes f_{j_1'}\otimes \cdots f_{j_m'}
\otimes \lambda \mu \otimes f_{l_1}\otimes \cdots \otimes
f_{l_n}\otimes \nu,$$ which becomes, after operating with
$U_{(k+m,s+t)(n,u)}$,
\begin{align*}
\sum_{j_1',\ldots, j_m'} u_{j_1',j_1}(s) \cdots u_{j_m',j_m}(s)
\sum_{l_1', \ldots l_n'} u_{l_1',l_1}(s+t) \cdots u_{l_n',l_n}(s+t)\\
 f_{i_1}\otimes \cdots \otimes f_{i_k}
\otimes f_{j_1'}\otimes \cdots f_{j_m'} \otimes f_{l_1'}\otimes
\cdots \otimes f_{l_n'}\otimes \lambda \mu \nu
\end{align*}
which is the same as above.

We now proceed to construct a product system representation that
will give rise to $\Psi$. We define
$$T_{(n,t)}(e_{i_1}\otimes \cdots e_{i_n} \otimes 1) = T_{i_1} \cdots T_{i_n}S_t .$$
The relation (\ref{unicomm}) is precisely what makes $T$ into a
product system representation (it is completely contractive
because $(T_i)_{i \in I\cup I'}$ is a row contraction). The last
thing to check is that for all $a \in B(H)$,
$$\tilde{T}_{(n,t)}(I_{X_{(n,t)}}\otimes a)\tilde{T}^*_{(n,t)} = \Psi_{(n,t)}(a).$$
But, after some identifications, $\tilde{S}_t = S_t$, and
$\tilde{T}$ is just the row contraction $(T_i)_{i\in I\cup I'}$,
so we are done.

We note that in this example too many ``miracles" have happened,
and we do not yet understand how what we have done here can be
generalized to other CP$_0$-semigroups over $\mathbb{N}\times
\Rp$.

\section{Appendix - examples of strongly commuting semigroups}
In this appendix we give some examples of strongly commuting
CP-semigroups. In special cases we are able to state a necessary
and sufficient condition for strong commutativity.

\subsection{Endomorphisms, automorphisms, and composition with automorphisms}
By Lemma 5.4 in \cite{S06}, there
are plenty of examples of CP maps $\Theta$,$\Phi$ that commute
strongly:
\begin{enumerate}
\item\label{it:end} If $\Theta$ and $\Phi$ are endomorphisms that commute then they commute strongly.
\item\label{it:aut} If $\Theta$ and $\Phi$ commute and either one of them is an automorphism then they commute strongly.
\item If $\alpha$ is a normal automorphism that commutes with $\Theta$, and $\Phi = \Theta \circ \alpha$, then $\Theta$ and $\Phi$ commute strongly.
\end{enumerate}

We note that item \ref{it:aut} does
not remain true if \emph{automorphism} is replaced by
\emph{endomorphism}. Because two CP-semigroups $\Theta$ and $\Phi$ commute strongly if
and only if for all $s,t\in \Rp$, $\Theta_s$ and $\Phi_t$ commute
strongly, it is immediate that:
\begin{enumerate}
\item\label{it:end1} If $\Theta$ and $\Phi$ are commuting E-semigroups then they commute strongly.
\item\label{it:aut1} If $\Theta$ and $\Phi$ commute and either one of them is an automorphism semigroup
then they commute strongly.
\item If $\alpha$ is a normal automorphism semigroup that commutes with
$\Theta$, and $\Phi_t = \Theta_t \circ \alpha_t$, then $\Theta$
and $\Phi$ commute strongly.
\end{enumerate}

At a first glance, item \ref{it:end1} might not seem very interesting in the
context of dilating CP-semigroups to enodmorphism semigroups. However, we find this this item \emph{very}
interesting, because one expcets a good dilation theorem not to complicate the situation in any sense.
For example, in Theorem \ref{thm:scudil}, in order to prove the existence of an E-dilation we have to assume that the CP-semigroups $\{R_t\}_{t\geq0}$ and $\{S_t\}_{t\geq0}$ are unital, but the E-dilation that we construct is also unital.
Another example, again from Theorem \ref{thm:scudil}: if the CP-semigroups act on a type $I$ factor, then so does the minimal E$_0$-dilation. The importance of item \ref{it:end1} is that it ensures that if $\{\alpha_t\}_{t\geq0}$ and $\{\beta_t\}_{t\geq0}$ are an E-dilation of $\{R_t\}_{t\geq0}$ and $\{S_t\}_{t\geq0}$, then $\alpha$ and $\beta$ commute strongly.

\subsection{Semigroups on $B(H)$}
It is a well known fact that if $\Theta$ and $\Phi$ are
CP-semigroups, then for each $t$ there are two
($\ell^2$-independent) row contractions $\{T_{i,t}\}_{i=1}^{m(t)}$ and
$\{S_{t,j}\}_{j=1}^{n(t)}$ ($m(t),n(t)$ may be equal to $\infty$) such that for all $a \in B(H)$
\be\label{eq:conj1}
\Theta_t(a) = \sum_{i}T_{t,i}aT_{t,i}^* ,
\ee
and
\be\label{eq:conj2}
\Phi_t(a) = \sum_{j}S_{t,j}aS_{t,j}^* \,.
\ee
We
shall call such semigroups \emph{conjugation semigroups}, as they
are given by conjugating an element with a row contraction. It now
follows from Proposition 5.8, \cite{S06}, that $\Theta$ and $\Phi$
commute strongly if and only if for all $(s,t)\in\Rpt$ there is an $m(t)n(t)\times m(t)n(t)$
unitary matrix
$$u(s,t) = \left(u(s,t)_{(i,j)}^{(k,l)}\right)_{(i,j),(k,l)}$$
such that for all $i,j$, \be\label{eq:SCunitary}T_{t,i}S_{s,j} =
\sum_{(k,l)}u(s,t)_{(i,j)}^{(k,l)}S_{s,l}T_{t,k} .\ee
As a simple example, if $\Phi$ and $\Psi$ are given by (\ref{eq:conj1}) and (\ref{eq:conj2}), and $S_{t,j}$ commutes with $T_{s,i}$ for all $s,t,i,j$, then $\Phi$ and $\Psi$ strongly commute.

\subsection{Semigroups on $B(H)$, $H$ finite dimensional}\label{subsec:H_finite}
If $H$ is a finite dimensional then any two commuting CP-semigroups on $B(H)$
commute strongly. This follows immediately from the following proposition.
\begin{proposition}\label{prop:SC_H_finite}
Let $\Phi$ and $\Psi$ be two commuting CP maps on $B(H)$, with $H$ a finite
dimensional Hilbert space. Then $\Phi$ and $\Psi$ strongly commute.
\end{proposition}
\begin{proof}
Assume that $\Phi$ is given by
$$\Phi(a) = \sum_{i=1}^m S_i a S_i^* $$
and that $\Psi$ is given by
$$\Psi(a) = \sum_{j=1}^n T_j a T_j^* ,$$
where $\{S_1, \ldots , S_m\}$ and $\{T_1, \ldots , T_n\}$ are row contractions
and $m,n \in \mathbb{N}$.
Because $\Phi$ and $\Psi$ commute, we have that
$$\sum_{i,j=1}^{mn} S_i T_j a T_j^* S_i^* = \sum_{i,j=1}^{mn} T_j S_i a S_i^* T_j^* $$
for all $a \in B(H)$. By the lemma on page 153 of \cite{G04} this implies that there exists an $mn \times mn$ unitary matrix $u$ such that
$$ S_i T_j = \sum_{(k,l)}u_{(i,j)}^{(k,l)}T_l S_k ,$$
and this means precisely that $\Phi$ and $\Psi$ strongly commute.
\end{proof}

We note here that the lemma cited above is stated in \cite{G04} for unital CP maps, but the proof works for the non-unital case as well. The reason that the assertion of the proposition fails for $B(H)$ with $H$ infinite dimensional is that in that case we may have $mn=\infty$, and the lemma is only true for a CP maps given by finite sums.

\subsection{Conjugation semigroups on general von Neumann algebras}
Let $\cM$ be a von Neumann algebra acting on a Hilbert space $H$.
We now show that if $\Theta$ and $\Phi$ are CP-semigroups on a von
Neumann algebra $\cM$ given as in (\ref{eq:conj1}) and
(\ref{eq:conj2}), where $T_{t,i},S_{t,j}$ are all in $\cM$, then a
sufficient condition for them to commute strongly with each other
is that there exists a unitary as in (\ref{eq:SCunitary}). To this
end, it is enough to show that if $\Theta$ and $\Phi$ are CP maps
given by
\bes
\Theta(a) = \sum_{i=1}^m T_{i}aT_{i}^* ,
\ees
and
\bes
\Phi(a) = \sum_{j=1}^n S_{j}aS_{j}^* ,
\ees
where $T_{i},S_{j}$ are
all in $\cM$, then a
sufficient condition for strong commutation is the existence of a
unitary matrix
$$u = \left(u_{(i,j)}^{(k,l)}\right)_{(i,j),(k,l)}$$
such that for all $i,j$, \bes T_{i}S_{j} =
\sum_{(k,l)}u_{(i,j)}^{(k,l)}S_{l}T_{k} .\ees Indeed, by
Proposition 5.6 of \cite{S06}, it is enough to show that there are
are two $\cM'$ correspondences $E$ and $F$, together with an
$\cM'$-correspondence isomorphism
$$t:E \otimes_{\cM'} F  \rightarrow F \otimes_{\cM'} E $$
and two c.c. representations $(\sigma,T)$ and $(\sigma,S)$ of $E$
and $F$, respectively, on $H$, such that:
\begin{enumerate}
\item\label{it:rep1} for all $a\in \cM$, $\tilde{T}(I_E \otimes
a)\tilde{T}^* = \Theta(a)$,
\item\label{it:rep2} for all $a\in \cM$, $\tilde{S}(I_F \otimes
a)\tilde{S}^* = \Phi(a)$,
\item\label{it:commute} $\tilde{T}(I_E \otimes \tilde{S}) = \tilde{S}(I_F \otimes
\tilde{T})\circ(t \otimes I_H)$.
\end{enumerate}
We construct these correspondences as follows. Let
$$E = \oplus_{i=1}^m \cM' \,\,\,\, {\rm and } \,\,\,\, F = \oplus_{j=1}^n \cM',$$
with the natural inner product and the natural actions of $\cM'$.
If we denote by $\{e_i \}_{i=1}^m$ and $\{f_j \}_{j=1}^n$ the
natural ``bases" of these spaces, then we can define
$$t(e_{i}\otimes f_{j}) =
\sum_{(k,l)}u_{(i,j)}^{(k,l)}f_{l}\otimes e_{k}.$$ We define
$\sigma$ to be the identity representation. Now $E \otimes_\sigma
H \cong \oplus_{i=1}^m H$, and $F \otimes_\sigma H \cong
\oplus_{j=1}^n H$, and on these spaces we define $\tilde{T}$ and
$\tilde{S}$ to be the row contractions given by $(T_1, \ldots ,
T_m)$ and $(S_1, \ldots, S_n)$. Some straightforward calculations
shows that items (\ref{it:rep1})-(\ref{it:commute}) are fulfilled.

\subsection{Semigroups on $\mathbb{C}^n$ or $\ell^\infty$}\label{subsec:SC_comm}
We close this paper with a more down-to-earth example of a
strongly commuting pair of CP$_0$-semigroups. Let $\cM = \mathbb{C}^n$ or
$\ell^\infty(\mathbb{N})$, considered as the algebra of diagonal
matrices acting on the Hilbert space $H = \mathbb{C}^n$ or
$\ell^2(\mathbb{N})$. In this context, a unital CP map is just a
stochastic matrix, that is, a matrix $P$ such that $p_{ij}\geq 0$
for all $i,j$ and such that for all $i$,
$$\sum_{j}p_{ij} = 1.$$
Indeed, it is straightforward to check that such a matrix gives rise to a
normal, unital, completely positive map. On the other hand, for all $i$, the composition
of a normal, unital, completely positive map with the normal state projecting onto the $i$th element
must be a normal state, so it has to be given by a nonnegative element in $\ell^1$ with norm $1$.

Given two such matrices $P$ and $Q$, we ask when do they strongly
commute. To answer this question, we first find orthonormal bases
for $\cM \otimes_P \cM \otimes_Q H$ and $\cM \otimes_Q \cM
\otimes_P H$. If $\{e_i\}$ is the vector with $1$ in the $i$th
place and $0$'s elsewhere, it is easy to see that the set $\{e_i
\otimes_P e_j \otimes_Q e_k\}_{i,j,k}$ spans $\cM \otimes_P \cM
\otimes_Q H$, and $\{e_i \otimes_Q e_j \otimes_P e_k\}_{i,j,k}$
spans $\cM \otimes_Q \cM \otimes_P H$. We compute
\begin{align*}
\langle e_i \otimes_P e_j \otimes_Q e_k, e_m \otimes_P e_p
\otimes_Q e_q \rangle &= \langle e_k, Q(e_j^*P(e_i^*
e_m)e_p)e_q\rangle \\
&=  \delta_{i,m} \delta_{j,p} \delta_{k,q}q_{kj}p_{ji}.
\end{align*}
Thus,
$$\{(q_{kj}p_{ji})^{-1/2} \cdot e_i \otimes_P e_j \otimes_Q e_k
: i,j,k \,\,{\rm such \,\, that}\,\, q_{kj}p_{ji} \neq 0\}$$ is an
orthonormal basis for $\cM \otimes_P \cM \otimes_Q H$, and
similarly for $\cM \otimes_Q \cM \otimes_P H$. If $u: \cM
\otimes_P \cM \otimes_Q H \rightarrow \cM \otimes_Q \cM \otimes_P
H$ is a unitary that makes $P$ and $Q$ commute strongly, then for
all $i,k$ we must have
$$u(e_i\otimes_P a \otimes_Q e_k) = (e_i\otimes 1 \otimes e_k)u(e_i\otimes_P a \otimes_Q e_k)
= e_i \otimes_Q b \otimes_P e_k ,$$ thus for all $i,j$, the spaces
$V_{i,j}:=\{e_i\otimes_P a \otimes_Q e_k : a \in \cM\}$ and $W_{i,j}:=\{e_i\otimes_Q
a \otimes_P e_k : a \in \cM\}$ bust be isomorphic. Thus, a
necessary condition for strong commutativity is that for all
$i,k$,
\be\label{eq:card}
|\{j :  q_{kj}p_{ji} \neq 0\} | = |\{j :
p_{kj}q_{ji} \neq 0\} | ,
\ee
where $| \cdot |$ denotes
cardinality. This condition is also sufficient, because we may define a
unitary between each pair $V_{i,j}$ and $W_{i,j}$, sending $e_i\otimes_P 1 \otimes_Q e_k$ to $e_i\otimes_Q 1 \otimes_P e_k$ and doing whatever on the complement. %
By the way, this example shows that when two CP maps commute strongly, there may be a great many unitaries
that ``implement" the strong commutation.

One can impose certain block structures on $P$ and $Q$ that will
guarantee that (\ref{eq:card}) is satisfied. Since we are in
particularly interested in semigroups, we shall be content with
the following observation. Let $P$ and $Q$ be two commuting, \emph{irreducible}, stochastic matrices.
Then $\cP_t := e^{-t}e^{tP}$ and $\cQ_t := e^{-t}e^{tQ}$ are two commuting, stochastic semigroups with strictly positive elements, and thus they commute strongly. For example, let
\[
P = \frac{1}{3}\left[
\begin{array}
[c]{ccc}%
1 & 1 & 1 \\
1 & 1 & 1 \\
1 & 1 & 1
\end{array}
\right]  \,\, , \,\,
Q = \left[
\begin{array}
[c]{ccc}%
{1}/{2} & 0 & {1}/{2} \\
{1}/{4} & {1}/{2} & {1}/{4} \\
{1}/{4} & {1}/{2} & {1}/{4}
\end{array}
\right]  .
\]
One may check that $P$ and $Q$ commute, but do not satisfy (\ref{eq:card}), hence they do not commute strongly. So we see that strong commutativity may fail even in the simplest cases.
However, $P$ and $Q$ are both irreducible, thus the semigroups they generate \emph{do} commute strongly.

\section{Acknowledgements}
The author is supported by the Jacobs School of Graduate Studies
and the Department of Mathematics at the Technion - I.I.T, and by
the Gutwirth Fellowship. This research is part of the author's
PhD. thesis, done under the supervision of Professor Baruch Solel.

\bibliographystyle{amsplain}

\end{document}